\documentclass[a4paper,11pt]{article}
\usepackage{latexsym,amsfonts,amsmath,mathrsfs}
\usepackage{graphicx}
\usepackage{color}
\usepackage[mathscr]{eucal}

\def\NN{\mathbb{N}}

\def\RR{\mathbb{R}}
\def\KK{\mathbb{K}}

\newcommand{\widebar}[1]{\overline{#1}}
\newcommand{\wb}[1]{\overline{#1}}

\newcommand{\wt}[1]{\widetilde{#1}}

\def\union{\mathop{\cup}}

\def\limsup{\mathop{\overline{\lim}}}
\def\sl{[\![}
\def\sr{]\!]}

\def\ds{\displaystyle}

\newcounter{hypot}

    \newenvironment{hypot}{\begin{list}
      {\hspace{\labelsep}\bfseries Assumption \Alph{hypot}.}
      {\leftmargin=0pt
       \labelwidth=0cm
       \refstepcounter{hypot}
       \def\makelabel##1{##1}}}{\end{list}}

\newcounter{assump}

\newenvironment{assump}{\begin{list}
      {\hspace{\labelsep}(\Alph{hypot}\arabic{assump})}
      {\leftmargin=1cm
       \labelwidth=1cm
       \usecounter{assump}
       \itshape }}{\end{list}}

\begin{document}

\bibliographystyle{plain}

\newtheorem{theorem}{Theorem}[section]
\newtheorem{proposition}[theorem]{Proposition}
\newtheorem{lemma}[theorem]{Lemma}
\newtheorem{corollary}[theorem]{Corollary}
\newtheorem{definition}[theorem]{Definition}
\newtheorem{remark}[theorem]{Remark}
\newtheorem{conjecture}[theorem]{Conjecture}
\newtheorem{assumption}[theorem]{Assumption}
\newtheorem{condition}[theorem]{Condition}
\newcommand{\defi}{\stackrel{\triangle}{=}}

\title{Impulsive control for continuous-time Markov Decision Processes}

\author{ \mbox{ }
F. Dufour
\thanks{Corresponding author.}
\\
\small Universit\'e Bordeaux I \\
\small IMB, Institut de Math\'ematiques de Bordeaux \\
\small INRIA Bordeaux Sud Ouest, Team: CQFD \\
\small 351 cours de la Liberation, 33405 Talence Cedex, France \\
\small e-mail : dufour@math.u-bordeaux1.fr
\and
A. B. Piunovskiy \\
\small Department of Mathematical Sciences \\
\small University of Liverpool  \\
\small L69 7ZL, Liverpool, United Kingdom \\
\small e-mail: piunov@liv.ac.uk
}

\maketitle

\begin{abstract}
The objective of this work is to study continuous-time Markov decision processes on a general Borel state space
with both impulsive and continuous controls for the infinite-time horizon discounted cost.
The continuous-time controlled process is shown to be non explosive under appropriate hypotheses.
The so-called Bellman equation associated to this control problem is studied.
Sufficient conditions ensuring the existence and the uniqueness of a bounded measurable solution to this optimality equation are provided.
Moreover, it is shown that the value function of the optimization problem under consideration satisfies this optimality equation. 
Sufficient conditions are also presented to ensure on one hand the existence of an optimal control strategy and on the other hand the existence of an $\varepsilon$-optimal control strategy.
The decomposition of the state space in two disjoint subsets is exhibited where roughly speaking, one should apply  a gradual action or an impulsive action correspondingly to get an optimal or $\varepsilon$-optimal strategy.
An interesting consequence of our previous results is as follows: the set of strategies that allow interventions at time $t=0$ and only immediately after natural jumps is a sufficient set for the control problem under consideration.
\end{abstract}
\begin{tabbing}
\small \hspace*{\parindent}  \= {\bf Keywords:}
Impulsive control, continuous control, continuous-time Markov decision process,\\
discounted cost.\\
\> {\bf AMS 2000 subject classification:} \= Primary 90C40; Secondary 60J25.
\end{tabbing}

\section{Introduction}
Continuous-Time Markov Decision Processes (CTMDP) form a general class of controlled stochastic processes. Those are suitable for formulating many optimization problems arising in such applications as engineering, computer science, telecommunications, finance, etc. The analysis of CTMDP started in the late 1950's and the early 1960's
with the pioneering works by Bellman, Blackwell, Howard, and Veinott, to name just a few authors; see, e.g., \cite{bellman57,howard60}. The analysis
has been mostly concentrated on control problems where the actions influence  the transition rate of the process continuously in time.
This is nowadays a very active area of research from the point of view of its theoretical foundations, as well as from the applications perspective, see e.g. the recent books and survey \cite{guo09,guo06,prieto12}.

Another class of models with impulsive actions, when the state of the process can be changed instantly, received very little attention.
The first attempt to study such problems is probably due to De Leve \cite{deleve64a,deleve64b}.
In the 1980's, a systematic study of impulsive control of continuous-time MDP, including a deterministic drift between the jumps,
has been conducted on one hand by Hordijk and Van Der Schouten and on the other hand by  Yushkevich.
Hordijk and Van Der Schouten \cite{hordijk83,hordijk84,hordijk85,schouten83} considered the case where only one impulsive action at a time is permitted. Given an observed history, the planned time moment for the next impulse was deterministic. In these papers, the optimization was performed within a special class of so called regular and conservative policies. One drawback of this approach is that the use of the dynamic programming method becomes problematic.
Yushkevich \cite{yushkevich83,yushkevich86,yushkevich87,yushkevich89} has introduced a new class of stochastic models, the so called T-processes where roughly speaking the processes are indexed by a parameter representing the natural current time and the number of the impulsive actions at that time moment. The introduction of this new family of processes was mainly motivated by the fact that it allows to consider models with multiple impulses at the same time moment.
For a general control model, Yushkevich typically studied the value functions of such control problems in terms of the related quasi-variational inequalities.
One should also mention another class of controlled models closely related to CTMDP and called piecewise deterministic Markov processes where impulsive control has also been considered.
Without attempting to present an exhaustive panorama, we advise the interested reader to consult the book \cite{davis93}, and references therein to get a rather complete view of this class of processes.

It is important to point out that impulsive control models are not mentioned in the recent monographs and surveys on CTMDP's \cite{guo09,guo06,prieto12}.
At the same time, they appear naturally in many real life situations. For example, in population dynamics, one can decide to remove an individual or immunize a susceptible in the epidemic framework. In the area of controlled queueing systems, one can decide to remove a job to give space for the customers with higher priority. In inventory or reliability theory, the decision maker can place an order or replace a device at any desired epoch.

The main difficulty of dealing with the impulsive control model is that impulsive actions give rise to a non standard path for the controlled process.
Indeed, the process may take several different values at the same time moment. This important property makes the classical theory of CTMDP's inapplicable.

\bigskip

The objective of the current paper is to develop a new approach to CTMDP's on a general Borel state space $\mathbf{X}$
with both impulsive and continuous controls.
In our framework, the continuous control influences the intensity of jumps $q$ at all times. This is in opposition to the impulsive control that intervenes by moving the process to a new point of the state space $\mathbf{X}$ at some moment specified by the controller.
In this context, continuous actions, also called gradual actions by Yushkevich (see for example \cite{yushkevich86}), take values in the space $\mathbf{A}^{g}$ and lead to \textit{natural} jumps, 
in opposition to an intervention of the controller on the process giving rise to an impulse.
In the latter case, at any time moment, one can apply an action from the set $\mathbf{A}^{i}$ of impulsive actions to change instantly
the state of the process according to a prescribed stochastic kernel $Q$ on $\mathbf{X}$ given $\mathbf{X}\times \mathbf{A}^{i}$. 
An intervention can lead not only to one single impulse but to any finite sequence of instantaneous impulses at the same time moment.
As a result, the controlled process can take several different values at the same time moment, the intervention epoch.
In the works of Hordijk and Van der Duyn Schouten \cite{hordijk83,hordijk84,hordijk85,schouten83}, only one impulsive action at a time was allowed.
As a consequence, the trajectory of the process was really a function of time, even if the intervention occurred immediately after natural jump.
In the works of Yushkevich \cite{yushkevich83,yushkevich86,yushkevich87,yushkevich89}, the time scale has been modified and split to make the trajectories as functions of time.
Therefore, a new theory of random processes had to be developed.
On the opposite, our aim is to use the standard theory of stochastic point processes \cite{bremaud81,jacod75,jacod79,last95}. In this context, it is necessary to extend the state space to take into account the fact that the controlled system may have several different values at the same time moment.
Our construction is based on a point process $(\Theta_{n},Y_{n})_{n\in \NN}$ where $\Theta_{n}$ represents the sojourn time between two consecutive epochs induced either by a natural jump or by an intervention. $Y_{n}$ is the new state vector of the form
\begin{eqnarray}
\label{defy}
(x_0,a_0,x_1,a_1,\ldots,x_k,a_k,x_{k+1},\Delta,\Delta,\ldots),
\end{eqnarray}
where $x_{0}$ corresponds to a possibly natural jump or to the value of the process just before the intervention. The pair $(a_{j},x_{j})$
(for $j\geq 1$) indicates that the impulsive action $a_{j}$ has been applied to the system, leading to a new location (jump) of the process denoted by $x_{j}$.
The special impulsive action $\Delta$ means that the impulses are over and the artificial state $\Delta$ means the same.
The space of all possible extended states as presented in (\ref{defy}) is denoted by $\mathbf{Y}$ (this set will be precisely defined in the next section).
The space of extended states resulting from interventions is denoted by $\mathbf{Y}^{*}=\mathbf{Y}\setminus\{(x_0,\Delta,\Delta,\ldots),~x_0\in X\}$.
Observe that $y=(x_0,\Delta,\Delta,\ldots)$ means no impulsive actions have been applied after a natural jump to state $x_0$.

We now present an informal description of the mechanism defining the controlled process $(\Theta_{n},Y_{n})_{n\in \NN}$.
In our framework, the interventions and gradual controls are determined through probability distributions on the appropriate spaces $\mathbf{Y}$ and
$\mathbf{A}^{g}$.
The initial time moment $0$ is very special. The initial state at the time moment just before $0$ is fixed and given by $Y_0=(x_0,\Delta,\Delta,\ldots)$
where $x_{0}$ is the initial location of the process. Moreover, the first sojourn time $\Theta_1$ equals zero.
Then, the controller chooses a probability measure on $\mathbf{Y}$ generating the random variable $Y_{1}$ which is
the next state immediately after time $0$.
 After this initial procedure, the controlled process can be constructed iteratively.
 Having observed the state $Y_{n}$, the controller chooses the action $u_{n}$ with the following components:
\begin{itemize}
\item a probability distribution on $\bar\RR^{*}_{+}$ generating the time of the (possible) next intervention which happens only in case no natural jumps occur earlier;
\item a stochastic kernel on $\mathbf{A}^{g}$ given $\RR_{+}$ describing the gradual control influencing the time of the next (possible) natural jump and its associated location;
\item an intervention immediately after the natural jump, in case it happens before the planned intervention, that is, a probability distribution on $\mathbf{Y}$;
\item a planned intervention, that is, a probability distribution on $\mathbf{Y}^{*}$. This last component is absent in the event that no interventions are allowed in the current state.
\end{itemize}
If the gradual action $a\in \mathbf{A}^{g}$ is applied at the state $x\in X$, then the cost rate is $C^g(x,a)$; any impulsive action $a\in \mathbf{A}^{i}$ results in the immediate cost $c^i(x,a)$. 
In the present paper, we consider the discounted model on the infinite time horizon. Note that an intervention occurs at one time moment with a fixed value of the discounting coefficient so that it corresponds to a discrete-time MDP with a total expected cost.

\bigskip

Our model is closely related to those studied by Hordijk, Van der Duyn Schouten and Yushkevich but presents important differences that we would like to emphasize. In particular,
in  \cite{yushkevich83,yushkevich86,yushkevich87,yushkevich89} only nonrandomized gradual controls were considered. Moreover, in \cite{hordijk83,hordijk84,hordijk85,schouten83,yushkevich83,yushkevich86,yushkevich87,yushkevich89} the authors consider the times of intervention as stopping times with respect to the filtration generated by the controlled process. In our context the times of intervention were specified through probability distributions depending on the history of the process.
In \cite{davis93} and the references therein, the control strategies were past-history independent, deterministic, and several impulses at the same time moment were forbidden.
Our framework is more general in the sense that we allow randomized policies.
Moreover, we allow instantaneous series of impulses which is not the case in \cite{davis93,hordijk83,hordijk84,hordijk85,schouten83}.
We would like to emphasize that \cite{yushkevich87} is the closest reference to our work because series of impulses is allowed.
The author studied the discounted cost control problem and showed that the value function is universally measurable and satisfies the Bellman equation.
Moreover,  the existence of an $\varepsilon$-optimal control strategy was proved.

\bigskip

When compared to the literature, our main contributions can be summarized as follows.
Our main objective in this paper is to study the Bellman equation associated with this control problem and to establish the existence of optimal and $\varepsilon$-optimal control strategies.
We first show that under some hypotheses the continuous-time controlled process is non explosive.
We provide sufficient conditions ensuring the existence and the uniqueness of a bounded measurable solution to the Bellman equation.
It is proved that this solution can be calculated by successive iterations of the associated Bellman operator.
Moreover, we show that the value function of our optimization problem satisfies this optimality equation. 
Two different sets of sufficient conditions are presented to ensure on one hand the existence of an optimal control strategy and on the other hand the existence of an $\varepsilon$-optimal control strategy.
An interesting consequence of our previous results is as follows: the set of strategies that allow intervention at time $t=0$ and only immediately after natural jumps is a sufficient set for the control problem under consideration.

To illustrate our theoretical results, we investigate the epidemic with carriers, such as typhoid. This model was suggested by Weiss in \cite{weiss65}
and was investigated by many authors \cite{booth86,piunovskiy04}. Similarly to \cite{piunovskiy04}, where the undiscounted model up to the end of the epidemic was considered, the optimal control strategy depends only on the number of the carriers and is of threshold type: immunize all the susceptibles as soon as the number of carriers exceeds a critical value $c^*$.

\bigskip

The rest of the paper is organized as follows. Section \ref{sec-process} is devoted to the construction of CTMDP's on a general Borel state space $\mathbf{X}$
with both impulsive and continuous controls.
Section \ref{sec-cost} introduces the infinite-horizon performance criterion and several different classes of admissible strategies. Several preliminary results are also formulated here.
The analysis of the Bellman equation and the existence of optimal and $\varepsilon$-optimal control strategies are discussed in section \ref{sec-main}.
Finally, Section \ref{sec-example} is devoted to the presentation of the example illustrating the results developed in the paper.
\section{The continuous-time Markov control process}
\label{sec-process}
The main goal of this section is to introduce the notations, as well as the parameters defining the model, and to present the construction of the controlled process.
In particular a measurable space $(\Omega,\mathcal{F})$ consisting of the canonical sample paths of the multivariate point process $(\Theta_{n},Y_{n})$ is introduced.
Having defined the class of admissible strategies, we show the existence of a probability measure $\mathbb{P}_{x_{0}}^{u}$ with respect to which the controlled process
$(\Theta_{n},Y_{n})$ has the required conditional distributions.

The following notations will be used in this paper:
$\NN$ is the set of natural numbers including $0$, $\NN^{*}=\NN-\{0\}$, $\RR$ denotes the set of real numbers, $\RR_{+}$ the set of non-negative real numbers,
$\RR_{+}^{*}=\RR_{+}-\{0\}$, $\widebar{\RR}_{+}=\RR_{+}\union \{+\infty\}$ and $\widebar{\RR}_{+}^*=\RR_{+}^*\union \{+\infty\}$.
For any $q\in \NN$, $\NN_{q}$ is the set $\{0,1,\ldots,q\}$ and for any $q\in \NN^{*}$, $\NN_{q}^{*}$ is the set $\{1,\ldots,q\}$.
The term \textit{measure} will always refer to a countably additive, ${\RR}_{+}$-valued set function.
Let $X$ be a Borel space and denote by $\mathcal{B}(X)$ its associated Borel $\sigma$-algebra.
For any set $A$, $I_{A}$ denotes the indicator function of the set $A$.
The set of measures defined on $(X,\mathcal{B}(X))$ is denoted by $\mathbb{M}(X)_{+}$, and
$\mathcal{P}(X)$ is the set of probability measures defined on $(X,\mathcal{B}(X))$, and
$\mathcal{P}(X|Y)$ is the set of stochastic kernels on $X$ given $Y$ where $Y$ denotes a Borel space.
For any point $x\in X$, $\delta_{x}$ denotes the Dirac measure defined by $\delta_{x}(\Gamma)=I_{\Gamma}(x)$ for any $\Gamma \in \mathcal{B}(X)$.
The set of bounded real-valued measurable functions defined on $X$ is denoted by $\mathbb{B}(X)$.
Finally, the infimum over an empty set is understood to be equal to $+\infty$.

\subsection{Parameters of the model}
We will deal with a control model defined through the following elements:
\begin{itemize}
\item $\mathbf{X}$ is the state space, assumed to be a Borel space (i.e., a measurable subset of a complete and separable metric space).
\item $\mathbf{A}$ is the action space, assumed to be also a Borel space. $\mathbf{A}^{i}\in\mathcal{B}(\mathbf{A})$
(respectively $\mathbf{A}^{g}\in \mathcal{B}(\mathbf{A})$) is the set of impulsive (respectively gradual) actions satisfying
$\mathbf{A}=\mathbf{A}^i\cup \mathbf{A}^g$ with $\mathbf{A}^i\cap \mathbf{A}^g=\emptyset$.
\item The set of feasible actions in state $x\in \mathbf{X}$ is $\mathbf{A}(x)$, which is a nonempty measurable subset of $\mathbf{A}$.
Admissible impulsive and gradual actions in the state $x\in \mathbf{X}$ are denoted by $\mathbf{A}^i(x)=\mathbf{A}(x)\cap \mathbf{A}^i$ and
$\mathbf{A}^g(x)=\mathbf{A}(x)\cap \mathbf{A}^g$.
It is supposed that $$\KK^g=\{(x,a)\in \mathbf{X}\times \mathbf{A}:a\in \mathbf{A}^{g}(x)\}\in\mathcal{B}(\mathbf{X}\times \mathbf{A}^g)$$
and this set contains the graph of a measurable function from $\mathbf{X}$ to $\mathbf{A}^g$ 
(necessarily $\mathbf{A}^g(x)\ne\emptyset$ for all $x\in \mathbf{X}$) and that 
$$\KK^i=\{(x,a)\in \mathbf{X}\times \mathbf{A}^{i}:a\in \mathbf{A}^i(x)\}\in\mathcal{B}(\mathbb{X}^{i}\times \mathbf{A}^i)$$
where $\mathbb{X}^{i}=\{x\in \mathbf{X}: \mathbf{A}^i(x)\ne\emptyset\}\in\mathcal{B}(\mathbf{X})$ and $\KK^i$ contains the graph of a measurable function from
$\mathbb{X}^{i}$ to $\mathbf{A}$.
\item The stochastic kernel $Q$ on $\mathbf{X}$ given $\KK^{i}$ describes the result of an impulsive action. 
In other words, if $x\in \mathbb{X}^{i}$ and an impulsive action $a\in \mathbf{A}^i(x)$ is applied then the state of the process changes instantly according to the stochastic kernel $Q$.
\item The signed kernel $q$ on $\mathbf{X}$ given $\KK^g$ is the intensity of jumps governing the dynamic of the process between interventions.
It satisfies $q(\mathbf{X}|x,a)=0$ and  $q(\mathbf{X}\setminus\{x\}|x,a)\ge 0$ for any $(x,a)\in\KK^g$.
\end{itemize}

In our model, an intervention consists only of a finite sequence of pairs of impulsive action and associated jump.
Actually, this finite sequence can be equivalently described by an infinite sequence of pairs of state and action, where the pairs are set to the fictitious action and state after a finite step.
As a result, an intervention is an element of the set
$$\mathbf{Y}=\bigcup_{k\in\NN}\mathbf{Y}_k \text{ with }
\mathbf{Y}_k=(\mathbf{X}\times \mathbf{A}^i)^k\times(\mathbf{X}\times\{\Delta\})\times(\{\Delta\}\times\{\Delta\})^\infty,$$
where $\Delta$ will play the role of the fictitious state and action.
The dynamic of such sequences is governed by the Markov Decision Process (MDP) $\mathcal{M}^{i}$ defined by
$$\mathcal{M}^{i}=\big(\mathbf{X}_\Delta,\mathbf{A}^i_\Delta,(\mathbf{A}^i_\Delta(x))_{x\in \mathbf{X}_\Delta},Q_{\Delta}\big)$$
where
\begin{itemize}
\item $\mathbf{X}_{\Delta}$, $\mathbf{A}^{i}_{\Delta}$ and $\big(\mathbf{A}^{i}_{\Delta}(x))_{x\in \mathbf{X}_{\Delta}}$
are the new state and actions spaces augmented by the fictitious state $\Delta$:
$\mathbf{X}_{\Delta}=\mathbf{X}\union\{\Delta\}$, $\mathbf{A}^{i}_{\Delta}=\mathbf{A}^{i}\union\{\Delta\}$ and
$\mathbf{A}^{i}_{\Delta}(x)=\mathbf{A}^{i}(x)\union\{\Delta\}$ for $x\in \mathbf{X}$ and $\mathbf{A}^{i}_{\Delta}(\Delta)=\{\Delta\}$.
\item $Q_{\Delta}(.|x,a)=Q(.|x,a)$ for any $(x,a)\in \KK^{i}$ and $Q_{\Delta}(\{\Delta\}|x,a)=1$ otherwise.
\end{itemize}
For the model $\mathcal{M}^{i}$, according to the Ionescu Tulcea's Theorem (see Proposition C.10 in \cite{hernandez96}),
 there exists a unique strategic measure $\beta^b(\cdot|x)$ on
$(\mathbf{X}_\Delta\times \mathbf{A}^i_\Delta)^\infty$ associated with the policy $b$ and the initial distribution $\delta_x$.
Here and below, we use the standard terminology for MDP: a policy is a sequence of past-dependent distributions on the action space; a Markov non-randomized policy is a sequence $(\varphi_j^i)_{j\in\NN}$ of $\mathbf{A}^i_\Delta$-valued mappings on $\mathbf{X}_\Delta$, and so on.
 Observe that $\beta^b$ is in fact a stochastic kernel on $(\mathbf{X}_\Delta\times \mathbf{A}^i_\Delta)^\infty$ given $\mathbf{X}$, see Proposition C.10 in \cite{hernandez96}.
Since we only consider intervention as an element of $\mathbf{Y}$, we restrict the admissible policies to the set $\Xi$ satisfying
$\beta^b(\mathbf{Y}|x)=1$ for $b\in \Xi$.
In fact, we consider \textit{randomized} interventions and consequently an intervention is an element of 
$$\mathcal{P}^\mathbf{Y}=\{\beta \in \mathcal{P}(\mathbf{Y}|\mathbf{X}): \beta(\cdot|\cdot)=\beta^b(\cdot|\cdot) \text{ for some } b\in\Xi\},$$
and
$$\mathcal{P}^\mathbf{Y}(x)=\{\rho \in \mathcal{P}(\mathbf{Y}): \rho(\cdot)=\beta^b(\cdot|x) \text{ for some } b\in\Xi\}$$
is the set of feasible interventions in state $x\in \mathbf{X}$.
Observe that if an intervention is chosen in $\mathbf{Y}_{0}$, it means actually that the controller has not intervened on the process through impulsive actions.
For technical reasons, it appears necessary to introduce the set $\mathbf{Y}^{*}$ of \textit{real} interventions given by
$$\mathbf{Y}^{*}=\bigcup_{k=1}^{\infty}\mathbf{Y}_k .$$ The associated sets of real \textit{randomized} interventions are defined by
$$\mathcal{P}^{
\mathbf{Y}^{*}}=\{\beta \in \mathcal{P}(\mathbf{Y}|\mathbf{X}): \beta(\cdot|\cdot)=\beta^b(\cdot|\cdot) \text{ for some } b\in\Xi
\text{ and }\beta^b(\mathbf{Y}^{*}|x)=1,\mbox{ for any }x\in \mathbb{X}^{i}\}$$
and
$$\mathcal{P}^{\mathbf{Y}^{*}}(x)=\{
\rho \in \mathcal{P}(\mathbf{Y}): \rho(\cdot)=\beta^b(\cdot|x) \text{ for some } b\in\Xi
\text{ and }\beta^b(\mathbf{Y}^{*}|x)=1\}$$
for $x\in \mathbf{X}$. Note that $\mathcal{P}^{\mathbf{Y}^{*}}(x)=\emptyset$ if $x\notin \mathbb{X}^{i}$.

Finally, we end this subsection by introducing a \textit{projection} mapping that will be used repeatedly in the paper.
If $y\in \mathbf{Y}$ then there exists a unique $k\in\NN$ such that $y\in \mathbf{Y}_k$. The $(k+1)$-th component of $y$ is of the form
$(x_{k+1},\Delta)$ with $x_{k+1}\in \mathbf{X}$.
The $\mathbf{X}$-valued mapping  $\bar x$ on $\mathbf{Y}$ is defined by
\begin{eqnarray}
\label{xbar}
\bar x(y)=x_{k+1}.
\end{eqnarray}

\subsection{Construction of the process}
Having introduced the parameters of the model, we are now in position to construct the Markov controlled process.
Let
$$\mathbf{Y}_\infty=\mathbf{Y}\cup\{y_{\infty}\}, \quad \mathbf{H}_{n}=\mathbf{Y}\times(\RR_{+}\times \mathbf{Y})^n, \quad \Omega_{n}=\mathbf{H}_{n}\times(\{\infty\}\times\{y_{\infty}\})^\infty,$$
for $n\in \NN$. The canonical space $\Omega$ is defined as
$$\Omega=\bigcup_{n=1}^\infty \Omega_{n}\bigcup \big( \mathbf{Y}\times(\RR_{+}\times \mathbf{Y})^\infty \big)$$
and is endowed with its Borel $\sigma$-algebra denoted by $\mathcal{F}$. For notational convenience, $\omega\in\Omega$ will be represented as
  $$\omega=(y_0,\theta_1,y_1,\theta_2,y_2,\ldots).$$
Here $y_0=(x_0,\Delta,\Delta,\ldots)$ is the initial state of the controlled point process $\xi$ with values in $\mathbf{Y}$, defined below;
$\theta_1=0$ and $y_1\in \mathbf{Y}$ is the result of the initial intervention.
The components $\theta_{n}>0$ for $n\geq 2$ mean the sojourn times; $y_{n}$ denotes the result of an intervention (if $y_{n}\in \mathbf{Y}^{*}$) or corresponds to a natural jump
(if $y_{n}\in \mathbf{Y}\setminus \mathbf{Y}^{*})$). In case $\theta_{n}<\infty$ and $\theta_{n+1}=\infty$, the trajectory has only $n$ jumps and we put $y_m=y_{\infty}$ (artificial point) for all $m\ge n+1$.

The path up to $n\in\NN$ is denoted by $$h_{n}=(y_0,\theta_1,y_1,\theta_2,y_2,\ldots \theta_{n},y_{n})\in \mathbf{H}_{n}.$$
For $n\in\NN$, introduce the mappings $Y_{n}:~\Omega\to \mathbf{Y}_\infty$ by $Y_{n}(\omega)=y_{n}$ and, for $n\in\NN^{*}$, the mappings
$\Theta_{n}:~\Omega\to \wb{\RR}_{+}$ by $\Theta_{n}(\omega)=\theta_{n}$.
The sequence $(T_{n})_{n\in\NN^{*}}$ of $\wb{\RR}_{+}$-valued mappings is defined on $\Omega$ by 
  $$T_{n}(\omega)=\sum_{i=1}^n\Theta_i(\omega)=\sum_{i=1}^n\theta_i$$
and $T_\infty(\omega)=\lim_{n\to\infty}T_{n}(\omega)$.
For notational convenience, we denote by
  $$H_{n}=(Y_0,\Theta_1,Y_1,\ldots,\Theta_{n},Y_{n})$$
the $n$-term history process taking values in $\mathbf{H}_{n}$ for $n\in\NN$.

The random measure $\mu$ associated with $(\Theta_{n},Y_{n})_{n\in \NN}$ is a measure defined on $\RR^{*}_{+}\times \mathbf{Y}$ by
  $$\mu(\omega;dt,dy)=\sum_{n\ge 2}I_{\{T_{n}(\omega)<\infty\}}\delta_{(T_{n}(\omega),Y_{n}(\omega))}(dt,dy).$$
For notational convenience the dependence on $\omega$ will be ignored and instead of $\mu(\omega;dt,dy)$ it will be written
$\mu(dt,dy)$.
The filtration $(\mathcal{F}_t)_{t\in\RR_{+}}$ on $(\Omega,\mathcal{F})$ is defined by
  $$\mathcal{F}_t=\sigma\{H_1\}\vee\sigma\{\mu(]0,s]\times B):~s\le t,B\in\mathcal{B}(\mathbf{Y})\}.$$
Finally, we define the controlled process $\big\{\xi_t\big\}_{t\in \RR_{+}}$:
  $$\xi_t(\omega)=\left\{
\begin{array}{ll}
Y_{n}(\omega), & \mbox{ if } T_{n}\le t<T_{n+1} \mbox{ for } n\in\NN^{*}; \\
y_{\infty}, & \mbox{ if } T_\infty\le t,
\end{array}\right.$$
and $\xi_{0-}(\omega)=Y_0=y_0$ with $y_0=(x_0,\Delta,\Delta,\ldots)$.
Obviously, the controlled process $(\xi_{t})_{t\in \RR_{+}}$ can be equivalently described by the sequence $(\Theta_{n},Y_{n})_{n\in \NN}$.

\subsection{Admisible strategies and conditional distribution of the controlled process}
An admissible control strategy is a sequence $u=(u_{n})_{n\in\NN}$ such that
$u_{0}\in \mathcal{P}^{\mathbf{Y}}(x_0)$ and, for any $n\in \NN^{*}$, $u_{n}$ is given by
$$u_{n}=\big( \psi_{n},\pi_{n},\gamma^0_{n},\gamma^1_{n} \big),$$
where
$\psi_{n}$ is a stochastic kernel on $\wb{\RR}^{*}_{+}$ given $\mathbf{H}_{n}$
satisfying $\psi_{n}(\cdot|h_{n})=\delta_{+\infty}(\cdot)$
for any $h_{n}=(y_0,\theta_1,\ldots \theta_{n},y_{n})\in \mathbf{H}_{n}$ with $\wb{x}(y_{n})\notin \mathbb{X}^{i}$,
$\pi_{n}$ is a stochastic kernel on $\mathbf{A}^{g}$ given $ \mathbf{H}_{n}\times \RR_{+}$
satisfying $\pi_{n}(\mathbf{A}^{g}(\wb{x}(y_{n}))|h_{n},t)=1$ for any $t\in \RR_{+}$ and
$h_{n}=(y_0,\theta_1,\ldots \theta_{n},y_{n})\in \mathbf{H}_{n}$,
$\gamma^0_{n}$ is a stochastic kernel on $\mathbf{Y}$ given $ \mathbf{H}_{n}\times \RR_+\times \mathbf{X}$ satisfying
$\gamma^0_{n}(\cdot|h_{n},t,\cdot) \in \mathcal{P}^{\mathbf{Y}}$ for any   $h_{n}\in \mathbf{H}_{n}$ and $t\in\RR_+$,
and $\gamma^1_{n}$  is a stochastic kernel on $\mathbf{Y}$  given $ \mathbf{H}_{n}$
satisfying $\gamma^1_{n}(\cdot|h_{n})\in  \mathcal{P}^{\mathbf{Y}^{*}}(\wb{x}(y_{n}))$
for any $h_{n}=(y_0,\theta_1,\ldots \theta_{n},y_{n})\in \mathbf{H}_{n}$
with $\wb{x}(y_{n})\in \mathbb{X}^{i}$; if $\wb{x}(y_n)\notin \mathbb{X}^i$ then $\gamma^1_n(\cdot|h_n)=\delta_{(\wb{x}(y_n),\Delta,\Delta,\ldots)}(\cdot)$.

\bigskip

The set of admissible control strategies is denoted by $\cal U$.
In what follows, we use notation $\gamma_{n}=(\gamma^0_{n},\gamma^1_{n})$.

\bigskip

Suppose a strategy $u=(u_{n})_{n\in\NN}\in \mathcal{U}$ is fixed with
$u_{n}=\big( \psi_{n},\pi_{n},\gamma^0_{n},\gamma^1_{n} \big)$ for $n\in \NN^{*}$.
We introduce the intensity of the natural jumps
\begin{eqnarray*}
\lambda_{n}(\Gamma_x,h_{n},t) & = & \int_{\mathbf{A}^{g}} \widebar{q}(\Gamma_x | \wb{x}(y_{n}),a) \pi_{n}(da | h_{n},t),
\end{eqnarray*}
where $\widebar{q}(\Gamma_x |x,a) = q(\Gamma_x\setminus\{x\}|x,a)$ for $(x,a)\in \mathbf{X}\times \mathbf{A}^{g}$, and
the rate of the natural jumps
\begin{eqnarray*}
\Lambda_{n}(\Gamma_{x},h_{n},t) & = & \int_{]0,t]} \lambda_{n}(\Gamma_{x},h_{n},s) ds 
\end{eqnarray*}
for any $n\in \NN^{*}$, $\Gamma_{x}\in \mathcal{B}(\mathbf{X})$ and $h_{n}=(y_0,\theta_1,y_1,\ldots,\theta_{n},y_{n})\in \mathbf{H}_{n}$.
Now, for any $n\in \NN^{*}$, the stochastic kernel $G_{n}$ on $\mathbf{Y}_{\infty}\times \wb{\RR}_{+}$ given $\mathbf{H}_{n}$ is defined by
\begin{eqnarray}
G_{n}(\{+\infty\}\times \{y_{\infty}\} | h_{n}) & = & \delta_{y_{n}} (\{y_{\infty}\}) + \delta_{y_{n}} (\mathbf{Y}) 
e^{-\Lambda_{n}(\mathbf{X},h_{n},+\infty)}\psi_{n}(\{+\infty\}|h_n)
\label{defG1}
\end{eqnarray}
and
\begin{align}
G_{n}(\Gamma_{\Theta} \times \Gamma_{y}| h_{n}) = \delta_{y_{n}} (\mathbf{Y}) \Big[ & 
\gamma_{n}^{1}(\Gamma_{y}| h_{n}) \int_{\Gamma_{\theta}} e^{-\Lambda_{n}(\mathbf{X},h_{n},t)} \psi_{n}(dt | h_{n}) \nonumber \\
& + \int_{\Gamma_{\theta}} \int_{\mathbf{X}} \psi_{n}([t,\infty] | h_{n}) \gamma_{n}^{0}(\Gamma_{y}| h_{n},t,x)  \lambda_{n}(dx,h_{n},t) e^{-\Lambda_{n}(\mathbf{X},h_{n},t)} dt \Big],
\label{defG2}
\end{align}
where $\Gamma_{y}\in \mathcal{B}(\mathbf{Y})$, $\Gamma_{\Theta}\in \mathcal{B}(\RR_{+})$
and $h_{n}=(y_0,\theta_1,y_1,\ldots,\theta_{n},y_{n})\in \mathbf{H}_{n}$.
Note that the kernel $\gamma^1_{n}$ does not appear in the formula for $G_{n}$ if $\wb{x}(y_{n})\notin \mathbb{X}^{i}$.

\bigskip

Consider an admissible strategy $u\in \mathcal{U}$ and an initial state $x_{0}\in\mathbf{X}$. From Theorem 3.6 in \cite{jacod75} (or Remark 3.43 in \cite{jacod79}), there exists a probability $\mathbb{P}^{u}_{x_{0}}$ on $(\Omega,\mathcal{F})$ such that
the restriction of $\mathbb{P}^{u}_{x_{0}}$ to $(\Omega,\mathcal{F}_{0})$ is given by
\begin{eqnarray}
\mathbb{P}^{u}_{x_{0}} \big( \{Y_{0}\}\times \{0\} \times \Gamma_y \times (\widebar{\RR}_{+}\times \mathbf{Y}_{\infty})^{\infty} \big) & = &
u_{0}(\Gamma_y|x_{0})
\label{Prob-init}
\end{eqnarray}
for any $\Gamma_y\in \mathcal{B}(\mathbf{Y})$ and the positive random measure $\nu$ defined on $\RR_{+}^{*}\times \mathbf{Y}$ by
\begin{eqnarray}
\nu(dt,dy)= \sum_{n\in \NN^{*}} \frac{G_{n}(dt-T_{n}, dy | H_{n})}{G_{n}([t,+\infty]\times \mathbf{Y}_{\infty} | H_{n})} I_{\{T_{n}< t \leq T_{n+1}\}}
\label{defnu}
\end{eqnarray}
is the predictable projection of $\mu$ with respect to $\mathbb{P}^{u}_{x_{0}}$.

\begin{remark}
\label{Cond-dist}
Observe that $\mathcal{F}_{T_{n}}$ is the $\sigma$-algebra generated by the random variable $H_{n}$ for $n\in \NN^{*}$.
The conditional distribution of $(Y_{n+1},\Theta_{n+1})$ given $\mathcal{F}_{T_{n}}$ under $\mathbb{P}^{u}_{x_{0}}$ is determined by 
$G_{n}(\cdot | H_{n})$ and the conditional survival function of $\Theta_{n+1}$ given $\mathcal{F}_{T_{n}}$ under $\mathbb{P}^{u}_{x_{0}}$
is given by  $G_{n}([t,+\infty]\times \mathbf{Y}_{\infty}| H_{n})$.
\end{remark}

\section{Optimization problem and preliminary results}
\label{sec-cost}
The objective of this section is to introduce the infinite-horizon performance criterion we are concerned with and several different classes of admissible strategies.
Some preliminary results are established. 
In particular, assuming the process is non explosive, a discounted version of the so-called Dynkin formula associated with the controlled process is derived (see Lemma \ref{lem1}).

\bigskip

The first result provides a decomposition of the predictable projection $\nu$ of the process in terms of two parts: one being related to the component
$(\gamma^{0}_{n})_{n\in \NN^{*}}$ of an admissible control strategy and the other to the component $(\gamma^{1}_{n})_{n\in \NN^{*}}$ 
\begin{lemma}
The predictable projection of the random measure $\mu$ is given by
\begin{eqnarray}
\nu=\nu_{0}+\nu_{1}
\label{nu=}
\end{eqnarray}
where for $\Gamma\in \mathcal{B}(\RR_{+}^{*})$, $\Gamma_{y}\in \mathcal{B}(\mathbf{Y})$
\begin{eqnarray}
\nu_{0}(\Gamma,\Gamma_{y}) & = & \int_{\Gamma}\int_{\mathbf{A}^{g}}\int_\mathbf{X} \gamma^{0}(\Gamma_{y}|x,s) \widebar{q}(dx | \wb{x}(\xi_{s-}),a) \pi(da|s) ds,
\label{nu0=} \\
\nu_{1}(\Gamma,\Gamma_{y}) & = & \sum_{n\in \NN^{*}} 
\gamma_{n}^{1}(\Gamma_{y}| H_{n}) \int_{\Gamma} I_{\{T_{n}< s \leq T_{n+1}\}}  \frac{\psi_{n}(ds-T_{n} | H_{n})}{\psi_{n}([s-T_{n},+\infty] | H_{n})},
\label{nu1=}
\end{eqnarray}
with
$$\gamma^{0}(dy|x,t)=\sum_{n\in \NN^{*}} I_{\{T_{n}< t \leq T_{n+1}\}}  \gamma^{0}_{n}(dy|H_{n},t-T_n,x),$$
and
$$\pi(da|t)=\sum_{n\in \NN^{*}} I_{\{T_{n}< t \leq T_{n+1}\}}  \pi_{n}(da|H_{n},t-T_{n})$$
for $t\in\RR_{+}$.
\end{lemma}
\textbf{Proof:}
First observe that by using the integration by parts formula, we obtain that
\begin{eqnarray*}
G_{n}([t,+\infty]\times \mathbf{Y}_{\infty}| h_{n}) & = & \delta_{y_{n}} (\{y_{\infty}\}) 
+ \delta_{y_{n}} (\mathbf{Y}) e^{-\Lambda_{n}(\mathbf{X},h_{n},t)} \psi_{n}([t,+\infty] | h_{n}).
\end{eqnarray*}
Now, recalling the definition of $\nu$ (see equation (\ref{defnu})) in terms of $G$ (see equation (\ref{defG2})), a straightforward calculation gives the result.
\hfill $\Box$

\bigskip

The cost rate $C^{g}$ associated with a gradual action is a real-valued mapping defined on $\mathbb{K}^{g}$.
The cost associated with an intervention  $y=(x_{0},a_{0},x_{1},a_{1},\ldots)\in \mathbf{Y}$ is given by 
\begin{equation}\label{ecd}
C^i(y)=\sum_{k\in \NN} c^i(x_k,a_k),
\end{equation}
where $c^{i}$ is a real-valued mapping defined on $\mathbf{X}_{\Delta}\times \mathbf{A}^{i}_{\Delta}$
satisfying $c^{i}(x,a)=0$ if $(x,a)\notin \mathbb{K}^{i}$.
For any $(x,a)\in \mathbb{K}^{i}$, $c^{i}(x,a)$ corresponds to the cost associated with a single jump at $x\in \mathbf{X}$  resulting from the impulsive action $a\in \mathbf{A}^{i}(x)$.
The cost associated with a randomized intervention $\beta\in \mathcal{P}^{\mathbf{Y}}(x)$ for $x\in \mathbf{X}$ is given by
$$\int_{\mathbf{Y}} C^i(y)\beta(dy|x).$$ 
Therefore, the infinite-horizon discounted performance criterion corresponding to an admissible control strategy $u\in{\cal U}$ is defined by
\begin{align}
\mathcal{V}(u,x_{0}) = \int_{\mathbf{Y}} C^{i}(y) & u_{0}(dy|x_{0}) + \mathbb{E}^{u}_{x_{0}} \Bigg[  \int_{0}^{+\infty} e^{-\eta s}  \int_{\mathbf{A}^{g}}  C^{g}(\wb{x}(\xi_{s-}),a) \pi(da |s) ds \Bigg] \nonumber \\
& + \mathbb{E}^{u}_{x_{0}} \Bigg[  \int_{]0,\infty[\times {\mathbf{Y}}} e^{-\eta s} C^{i}(y) \mu(ds,dy) \Bigg].
\label{Def-cost}
\end{align}
In the previous expression, $\eta>0$ is the discount factor, $\mathcal{V}(u,x_{0})$ is understood to be equal to $+\infty$ if the integrals of both the positive and negative parts of the integrand are infinite.
Note that, for any control strategy $u\in{\cal U}$, the function ${\cal V}(u,\cdot)$ is measurable.
The optimization problem under consideration is to minimize $\mathcal{V}(u,x_{0})$ within the class of admissible strategies $u\in \mathcal{U}$ where $x_{0}$ is the initial state.
A control strategy $u\in{\cal U}$ is called
\begin{itemize}
\item {\it non-randomized stationary}, if $\psi_n(\cdot|h_n)=\delta_{\psi^s(\wb{x}(y_n))}(\cdot)$, $\pi_n(\cdot|h_n,t)=\delta_{\varphi^s(\wb{x}(y_n))}(\cdot)$, $\gamma^0_n(\cdot|h_n,t,\cdot)=\beta^{b_0}(\cdot|\cdot)$ and $\gamma^1_n(\cdot|h_n)=\beta^{b_1}(\cdot|\wb{x}(y_n))$, where $\psi^s$ 
(respectivelt, $\varphi^s$) is a measurable map from $\mathbf{X}$ to $\wb{\RR}^{*}_{+}$ (respectively, from $\mathbf{X}$ to $\mathbf{A}^g$) and
$b_0,b_1$ are non-randomized stationary policies in ${\cal M}^i$.
\item {\it non-randomized almost stationary}, in case $b_0$ and $b_1$ in the above definition are Markov non-randomized policies.
\item {\it uniformly or persistently} optimal (respectively, $\varepsilon$-optimal for $\varepsilon>0$), if $\ds \mathcal{V}(u,x_{0})=\inf_{v\in{\cal U}} \mathcal{V}(v,x_{0})$ (respectively, $\mathcal{V}(u,x_{0})\le \mathcal{V}(v,x_{0})+\varepsilon$ for any $v\in{\cal U}$) simultaneously for all $x_0\in\mathbf{X}$ and hence for any initial distribution.
\end{itemize}

\bigskip

The following lemma provides a discounted version of the so-called Dynkin formula associated with the controlled process $(\xi_{t})_{t\in \RR_{+}}$
\begin{lemma}\label{lem1}
Suppose a strategy $u=(u_{n})_{n\in\NN}\in \mathcal{U}$ is fixed with
$u_{n}=\big( \psi_{n},\pi_{n},\gamma^0_{n},\gamma^1_{n} \big)$ for $n\in \NN^{*}$
satisfying $\mathbb{P}^{u}_{x_{0}}(T_{\infty}=+\infty)=1$.
Let $W$ be a bounded real-valued measurable function defined on $\mathbf{X}$ and $\eta>0$ be a discount factor.
Then
\begin{align}
& \mathbb{E}^{u}_{x_{0}}  \big[e^{-\eta t} W(\wb{x}(\xi_{t}))\big] =  \int_{\mathbf{Y}} W(\widebar{x}(y)) u_{0}(dy|x_{0})
+ \mathbb{E}^{u}_{x_{0}} \Bigg[
\int_{0}^{t} \int_{\mathbf{A}^{g}}
e^{-\eta s} \Big[ -\eta W(\wb{x}(\xi_{s})) \nonumber \\
& +\int_{\mathbf{X}} \int_{\mathbf{Y}} W(\widebar{x}(y)) \gamma^{0}(dy | x,s) \widebar{q}(dx| \wb{x}(\xi_{s}),a)
-W(\wb{x}(\xi_{s})) \widebar{q}(\mathbf{X}| \wb{x}(\xi_{s}),a)    \Big] \pi(da |s) ds \Bigg] \nonumber \\
& + \mathbb{E}^{u}_{x_{0}} \Bigg[ \sum_{n\in \NN^{*}}\int_{]T_{n}\wedge t, T_{n+1}\wedge t]}\int_{\mathbf{Y}} e^{-\eta s} \Big[ W(\widebar{x}(y)) - W(\wb{x}(\xi_{s-})) \Big]  
\gamma^{1}_{n}(dy |H_{n}) \frac{\psi_{n}(ds-T_{n} | H_{n})}{\psi_{n}([s-T_{n},+\infty] | H_{n})} \Bigg]
\end{align}
\end{lemma}
\textbf{Proof:}
By using the fact that $\mathbb{P}^{u}_{x_{0}}(T_{\infty}=+\infty)=1$ and the product formula for functions of bounded variation
(see for example Theorem A.4.6 in \cite{last95}) we have that
\begin{eqnarray*}
e^{-\eta t} W(\wb{x}(\xi_{t})) & = & W(\wb{x}(y_{1})) - \int_{0}^{t} \eta e^{-\eta s}  W(\wb{x}(\xi_{s})) ds \\
&&+ \int_{]0,t]\times \mathbf{Y}} e^{-\eta s} \Big[ W(\wb{x}(z)) - W(\wb{x}(\xi_{s-})) \Big] \mu(ds,dz) .
\end{eqnarray*}
Consequently, by using the fact that $\nu$ is the predictable projection of $\mu$ and $W$ is bounded, it yields
\begin{eqnarray*}
\mathbb{E}^{u}_{x_{0}} \big[ e^{-\eta t} W(\wb{x}(\xi_{t})) \big] & = & \mathbb{E}^{u}_{x_{0}} \big[ W(\wb{x}(y_{1}))  \big]
- \mathbb{E}^{u}_{x_{0}} \Bigg[
\int_{0}^{t} \eta e^{-\eta s}  W(\wb{x}(\xi_{s})) ds \Bigg] \\
& & + \mathbb{E}^{u}_{x_{0}} \Bigg[ \int_{]0,t]\times \mathbf{Y}} e^{-\eta s} \Big[ W(\wb{x}(z)) - W(\wb{x}(\xi_{s-})) \Big] \nu(ds,dz) 
\Bigg].
\end{eqnarray*}
Combining the previous equation with (\ref{nu=}), (\ref{nu0=}) and (\ref{nu1=}), we obtain
\begin{align*}
& \mathbb{E}^{u}_{x_{0}} \big[ e^{-\eta t} W(\wb{x}(\xi_{t})) \big] = \mathbb{E}^{u}_{x_{0}} \big[ W(\wb{x}(y_{1}))  \big]
- \mathbb{E}^{u}_{x_{0}} \Bigg[
\int_{0}^{t} \eta e^{-\eta s}  W(\wb{x}(\xi_{s})) ds \Bigg] \\
& + \mathbb{E}^{u}_{x_{0}} \Bigg[ \int_{]0,t]}  \int_{{\mathbf{Y}}} e^{-\eta s} \Big[ W(\widebar{x}(y)) - W(\wb{x}(\xi_{s-})) \Big]  
\gamma^{0}(dy|x,s) \int_{\mathbf{A}^{g}}\int_\mathbf{X} \widebar{q}(dx | \wb{x}(\xi_{s-}),a) \pi(da|s) ds \Bigg]\\
& + \mathbb{E}^{u}_{x_{0}} \Bigg[ \sum_{n\in \NN^{*}}\int_{]T_{n}\wedge t, T_{n+1}\wedge t] }\int_{{\mathbf{Y}}} e^{-\eta s} \Big[ W(\widebar{x}(y)) - W(\wb{x}(\xi_{s-})) \Big]  
\gamma^{1}_{n}(dy | H_{n}) \frac{\psi_{n}(ds-T_{n} | H_{n})}{\psi_{n}([s-T_{n},+\infty] | H_{n})} \Bigg].
\end{align*}
Now, from equation (\ref{Prob-init}), it follows that
\begin{eqnarray*}
\mathbb{E}^{u}_{x_{0}} \big[ W(\wb{x}(y_{1}))  \big] & = & \int_{\mathbf{Y}} W(\widebar{x}(y)) u_{0}(dy|x_{0})
\end{eqnarray*}
showing the result.
\hfill $\Box$

\section{Main results}
\label{sec-main}
This section is devoted to the analysis of the so-called Bellman equation associated with the control problem described in the previous section
and to the existence of optimal and $\varepsilon$-optimal control strategies.
The first result (see Proposition \ref{explo}) ensures that the continuous-time controlled process is non explosive under some hypotheses.
Then we provide two different sets of conditions (see Assumptions \ref{HC1} and \ref{HC1}) ensuring the existence of a bounded measurable solution to the Bellman equation.
More precisely, it is proved in Propositions \ref{existence-Bell} and \ref{existence-Bell2} that this solution can be calculated by the successive iteration of the associated Bellman operator, leading either to an upper semicontinuous or to a lower semicontinuous solution.
Moreover, we show in Theorem \ref{th1} and Corollary \ref{corol1}, on one hand, the existence of an optimal control strategy and, on the other hand, the existence of an $\varepsilon$-optimal control strategy. We also prove that the value function of the optimization problem under consideration satisfies this optimality equation and, as a consequence, the bounded solution of the Bellman equation is unique.
We exhibit the decomposition of the state space in two disjoint subsets $\mathbf{X}^i$ and $\mathbf{X}^g$ where, roughly speaking, one should apply a gradual action if the current state is in $\mathbf{X}^g$, and an impulsive action if the current state is in $\mathbf{X}^i$, to get an optimal or an $\varepsilon$-optimal strategy,
depending on the assumptions under consideration (see Remark \ref{decomp}).
Another important and interesting consequence of our previous results is as follows: the set of strategies that allow intervention at time $t=0$ and only immediately after natural jumps is a sufficient set for the control problem under study. (See Theorem \ref{th1} and Corollary \ref{corol1}.)

\bigskip

The Bellman equation reads as follows:
\begin{align}
\inf_{a\in \mathbf{A}^{g}(x)} & \Big\{ -\eta V(x) + \int_{\mathbf{X}} V(\tilde x) \widebar{q}(d\tilde x | x,a) -V(x) \widebar q (\mathbf{X}|x,a) + C^{g}(x,a) \Big\} \nonumber \\
& \wedge \inf_{a\in \mathbf{A}^{i}(x)}  \Big\{ -V(x) + \int_{\mathbf{X}} V(\tilde x) Q(d\tilde x | x,a) +c^{i}(x,a) \Big\} = 0
\label{Bell}
\end{align}
for any $x\in \mathbf{X}$.
If $V$ is a solution to the equation (\ref{Bell}), we introduce the following subsets of $\mathbf{X}$:
\begin{eqnarray}
\mathbf{X}^{g} & = & \Big\{ x\in \mathbf{X} : \eta V(x) =
\inf_{a\in \mathbf{A}^{g}(x)} \Big\{ \int_{\mathbf{X}} V(\tilde x) \widebar{q}(d\tilde x | x,a) -V(x) \widebar q (\mathbf{X}|x,a) + C^{g}(x,a) \Big\} \Big\},
\end{eqnarray}
and 
\begin{eqnarray}
\mathbf{X}^{i} & = &  \mathbf{X}\setminus\mathbf{X}^{g} \subset \Big\{ x\in \mathbf{X} : V(x) =
\inf_{a\in \mathbf{A}^{i}(x)} \Big\{\int_{\mathbf{X}} V(\tilde x) Q(d\tilde x | x,a) +c^{i}(x,a) \Big\} \Big\},
\end{eqnarray}
These sets will be used to construct an optimal or an $\varepsilon$-optimal strategy in the proof of Theorem \ref{th1}.
Below, we provide conditions under which there exists a measurable bounded solution to the Bellman equation. Those conditions also guarantee that the sets
$\mathbf{X}^{g}$ and $\mathbf{X}^{i}$ are measurable.

\begin{hypot}
\item\label{HA} There exists a constant $K\in \RR$ such that for any $x\in \mathbf{X}$ and $a^g\in \mathbf{A}^{g}(x)$ and
$a^i\in \mathbf{A}^{i}(x)$:
\begin{assump}
\item \label{HA1} $\widebar{q}(\mathbf{X}|x,a^g)\leq K$. 
\item \label{HA2} $\big| C^{g}(x,a^g) \big| \leq K$.
\item \label{HA4} $c^{i}(x,a^i)\geq 0$.
\end{assump}
\end{hypot}

\begin{hypot}
\item\label{HB}
There exists a constant  $\underline{c}>0$ such that $c^{i}(x,a)\geq \underline{c}$ for any $(x,a)\in \KK^{i}$.
\end{hypot}

The following proposition gives sufficient condition for non explosion.
\begin{proposition}
\label{explo}
Suppose that Assumptions \ref{HA} and \ref{HB} hold. If $u\in \mathcal{U}$ satisfies $\mathcal{V}(u,x_{0})<\infty$ then 
$\mathbb{P}^{u}_{x_{0}}(T_{\infty}<\infty)=0$.
\end{proposition}
\textbf{Proof:}
From Assumption \ref{HA} and the definition of the cost (\ref{Def-cost}), we have that
\begin{eqnarray*}
\mathcal{V}(u,x_{0}) & \geq & -\frac{K}{\eta} + \mathbb{E}^{u}_{x_{0}} \Bigg[  \int_{]0,\infty[\times \mathbf{Y}} e^{-\eta s} C^{i}(y) \mu(ds,dy) \Bigg]
\nonumber \\
& \geq & -\frac{K}{\eta} + \mathbb{E}^{u}_{x_{0}} \Bigg[ \sum_{n\in \NN^{*}}\int_{]T_{n}, T_{n+1}] \times \mathbf{Y}} e^{-\eta s}
C^{i}(y) \gamma^{1}_{n}(dy |H_{n}) \frac{\psi_{n}(ds-T_{n} | H_{n})}{\psi_{n}([s-T_{n},+\infty] | H_{n})} \Bigg].
\end{eqnarray*}
Now, observe that if $\wb{x}(Y_n)\notin \mathbb{X}^{i}$ then the measure 
$e^{-\eta s} \gamma^{1}_{n}(dy |H_{n}) \frac{\psi_{n}(ds-T_{n} | H_{n})}{\psi_{n}([s-T_{n},+\infty] | H_{n})}$ is zero on the set $]T_{n}, T_{n+1}] \times \mathbf{Y}$
and if $\wb{x}(Y_n)\in \mathbb{X}^{i}$ then $\gamma^1_n(\cdot|H_n)\in{\cal P}^{\mathbf{Y}^{*}}(\wb{x}(Y_n))$,
and that $C^i(y)\ge\underline{c}$ for any $y\in \mathbf{Y}^{*}$ by Assumption \ref{HB}. Consequently,
\begin{eqnarray}
\mathcal{V}(u,x_{0})  
& \geq & -\frac{K}{\eta} + \underline{c} \: \mathbb{E}^{u}_{x_{0}} \Bigg[ \sum_{n\in \NN^{*}}\int_{]T_{n}, T_{n+1}] \times \mathbf{Y}} e^{-\eta s}
\gamma^{1}_{n}(dy |H_{n}) \frac{\psi_{n}(ds-T_{n} | H_{n})}{\psi_{n}([s-T_{n},+\infty] | H_{n})} \Bigg].
\label{explo-1}
\end{eqnarray}

Moreover, from Assumption \ref{HA1} we get that
\begin{eqnarray}
\mathbb{E}^{u}_{x_{0}} \Bigg[ \int_{0}^{+\infty} \int_{\mathbf{A}^{g}}
e^{-\eta s} \int_{\mathbf{X}} \int_{\mathbf{Y}}   \gamma^{0}(dy | x,s) \widebar{q}(dx| \wb{x}(\xi_{s}),a) \pi(da |s) ds \Bigg] & \leq & \frac{K}{\eta}.
\label{explo-2}
\end{eqnarray}
Combining equations (\ref{explo-1}) and (\ref{explo-2}), we have that
\begin{eqnarray}
\mathbb{E}^{u}_{x_{0}} \Bigg[ \int_{0}^{+\infty} \int_{\mathbf{Y}} e^{-\eta s} \mu(ds,dz) \Bigg] = \mathbb{E}^{u}_{x_{0}} \Bigg[ \int_{0}^{+\infty} \int_{\mathbf{Y}}  e^{-\eta s} \nu(ds,dz) \Bigg] & \leq & \frac{1}{\underline{c}} \Big[\mathcal{V}(u,x_{0})+\frac{K}{\eta} \Big]+\frac{K}{\eta}.
\label{explo-3}
\end{eqnarray}
However, if $\mathbb{P}^{u}_{x_{0}}(T_{\infty}<\infty)>0$ then
\begin{eqnarray}
\mathbb{E}^{u}_{x_{0}} \Bigg[ \int_{0}^{+\infty} \int_{\mathbf{Y}} e^{-\eta s} \mu(ds,dz) \Bigg] & \geq &
\mathbb{E}^{u}_{x_{0}} \Big[ e^{-\eta T_{\infty}} \mu(\RR^*_+,\mathbf{Y}) I_{\{T_{\infty}<\infty\}} \Big] = +\infty.
\label{explo-4}
\end{eqnarray}
From equations (\ref{explo-3}) and (\ref{explo-4}), it follows that if $u\in \mathcal{U}$ satisfies $\mathcal{V}(u,x_{0})<\infty$ then
$\mathbb{P}^{u}_{x_{0}}(T_{\infty}<\infty)=0$, showing the result.
\hfill $\Box$

\bigskip

In Assumption \ref{HC} below, we assume that metrizable topologies in the spaces $\mathbf{X}$ and $\mathbf{A}$ are fixed.
\begin{hypot}
\item\label{HC} \mbox{ }
\begin{assump}
\item \label{HC1} The sets $\mathbb{K}^{g}$ and $\mathbb{K}^{i}$ are open in $\mathbf{X}\times\mathbf{A}^g$ and $\mathbb{X}^{i}\times\mathbf{A}^i$ correspondingly.
For any continuous bounded function $F$ on $\mathbf{X}$, the functions $\displaystyle \int_{\mathbf{X}}F(z)\widebar q(dz|x,a)$ and
$\displaystyle \int_{\mathbf{X}}F(z)Q(dz|x,a)$ are continuous on $\mathbb{K}^{g}$ and $\mathbb{K}^{i}$ correspondingly.
The functions $C^g$ and $c^i$ are upper semicontinuous on $\mathbb{K}^{g}$ and $\mathbb{K}^{i}$ correspondingly.
\item \label{HC2} The sets $\mathbf{A}^g$ and $\mathbf{A}^i$ are compact and the sets $\mathbb{K}^{g}$ and $\mathbb{K}^{i}$ are closed in $\mathbf{X}\times\mathbf{A}^g$ and $\mathbb{X}^{i}\times\mathbf{A}^i$ correspondingly.
For any continuous bounded function $F$ on $\bf X$, the functions $\displaystyle \int_{\mathbf{X}}F(z)\widebar q(dz|x,a)$ and $\displaystyle \int_\mathbf{X}F(z)Q(dz|x,a)$
are continuous on $\mathbb{K}^{g}$ and $\mathbb{K}^{i}$ correspondingly.
The functions $C^g$ and $c^i$ are lower semicontinuous on $\mathbb{K}^{g}$ and $\mathbb{K}^{i}$ correspondingly.
\end{assump}
\end{hypot}

\bigskip

Introduce the stochastic kernel $\widetilde{P}$ on $\mathbf{X}$ given $\mathbb{K}^{g}$
$$\widetilde{P}(\Gamma|x,a)=\frac{1}{K}\Big[\widebar{q}(\Gamma|x,a)+\delta_{x} (\Gamma) \big[K-\widebar{q}(\mathbf{X}|x,a)\big]\Big]$$
for any $\Gamma\in \mathcal{B}(\mathbf{X})$ and $(x,a)\in\mathbb{K}^{g}$
and consider the mapping $\mathfrak{B}$ defined on $\mathbb{B}(\mathbf{X})$ by
\begin{align}
\mathfrak{B}F(x)=\inf_{a\in \mathbf{A}^{g}(x)} & \Big\{ \frac{K}{K+\eta}\int_{\mathbf{X}} F(\widetilde x) \widetilde{P}(d\widetilde x |x,a)+ \frac{1}{K+\eta} C^{g}(x,a) \Big\} \nonumber \\
& \wedge \inf_{a\in \mathbf{A}^{i}(x)}  \Big\{ \int_{\mathbf{X}} F(\widetilde x) Q(d\widetilde x | x,a) + c^{i}(x,a) \Big\}
\label{Bell-a}
\end{align}
for any $F\in \mathbb{B}(\mathbf{X})$. The mapping $\mathfrak{B}$ will be called the \textit{Bellman operator} for further references.

\bigskip

The next two propositions ensure, under two different sets of conditions, the existence of an upper semicontinuous or a lower semicontinuous solution of the Bellman equation,
the measurability of the corresponding sets $\mathbf{X}^g$ and $\mathbf{X}^i$ and the existence of Borel-measurable mappings $\varphi^i:\mathbf{X}^i\to\mathbf{A}^i$ 
and $\varphi^g:\mathbf{X}^g\to\mathbf{A}^g$ that will be used to construct optimal strategies.
\begin{proposition}
\label{existence-Bell}  Suppose Assumptions \ref{HA} and \ref{HC1} hold.
Then the decreasing sequence of functions 
$(V_{i})_{i\in \NN}$ defined iteratively by $V_{i+1}=\mathfrak{B}V_{i}$ with $V_{0}=\frac{K}{\eta}$
belongs to $\mathbb{B}(\mathbf{X})$ and converges to a bounded upper semicontinuous function $V$ on $\mathbf{X}$ satisfying the Bellman equation (\ref{Bell}).
Moreover, the corresponding sets $\mathbf{X}^g$ and $\mathbf{X}^i$ are measurable
and, for any $\varepsilon>0$, there exist Borel-measurable mappings $\varphi^i:\mathbf{X}^i\to\mathbf{A}^i$ and $\varphi^g:\mathbf{X}^g\to\mathbf{A}^g$, such that
\begin{equation}\label{qu1}
\varphi^i(z)\in\left\{a\in\mathbf{A}^i(z):~\int_\mathbf{X}V(\widetilde x)Q(d\widetilde x|z,a)+c^i(z,a)\le V(z)+\varepsilon\right\},
\end{equation}
for any $z\in\mathbf{X}^i$ and
\begin{equation}\label{qu2}
\varphi^g(z)\in\left\{a\in\mathbf{A}^g(z):~\int_\mathbf{X}V(\widetilde x)\widebar q(d\widetilde x|z,a)-V(z)\widebar q(\mathbf{X}|z,a)+C^g(z,a)\le \eta V(z)+\varepsilon\right\},
\end{equation}
for any $z\in\mathbf{X}^g$.
\end{proposition}
\textbf{Proof:}
By using simple algebraic manipulations and Assumptions \ref{HA1}-\ref{HA2}, it is easy to show that
$V\in \mathbb{B}(\mathbf{X})$ is a solution of the Bellman equation (\ref{Bell}) if and only if  $V\in \mathbb{B}(\mathbf{X})$ and satisfies $V=\mathfrak{B}V$.
Let us denote by $\mathbb{U}(\mathbf{X})$ the set of upper semicontinuous functions defined on $\mathbf{X}$.
Clearly, from Proposition 7.34 in \cite{bertsekas78} and Assumption \ref{HC1}, the operator $\mathfrak{B}$ maps $\mathbb{U}(\mathbf{X})$ into $\mathbb{U}(\mathbf{X})$.
Consider the sequence $(V_{i})_{i\in \NN}$  defined by $V_{i+1}=\mathfrak{B}V_{i}$ with $V_{0}=\frac{K}{\eta}$.
We will show that $V_i\in\mathbb{B}(\mathbf{X})$ for any $i\in\NN$.
By definition of $\mathfrak{B}$ and Assumptions \ref{HA1}-\ref{HA2}, we have
\begin{eqnarray*}
V_{1}(x) & \leq & \inf_{a\in \mathbf{A}^{g}(x)} \Big\{ \frac{K}{K+\eta}\int_{\mathbf{X}} V_{0}(\widetilde x) \widetilde{P}(d\widetilde x |x,a)+ \frac{1}{K+\eta} C^{g}(x,a) \Big\} \\
& \leq & \frac{K}{K+\eta} \frac{K}{\eta} + \frac{K}{K+\eta}=\frac{K}{\eta}=V_{0}(x).
\end{eqnarray*}
From the previous inequality and since the operator $\mathfrak{B}$ is monotone, it can be easily shown by induction that the sequence $(V_{i})_{i\in \NN}$
belongs to $\mathbb{U}(\mathbf{X})$ and satisfies 
\begin{eqnarray}
V_{i+1} =\mathfrak{B} V_{i} \leq V_{i},
\label{ineqVi}
\end{eqnarray}
for any $i\in \NN$.
Moreover, we have that  $\sup_{x\in \mathbf{X}}|V_{i}(x)|\leq \frac{K}{\eta}$.
Indeed, from equation (\ref{ineqVi}), it follows easily that $V_{i}(x)\leq \frac{K}{\eta}$.
Let us show by induction that $V_{i}(x)\geq -\frac{K}{\eta}$.
Clearly, we have $V_{0}(x)\geq -\frac{K}{\eta}$.
Assume that $V_{i}(x)\geq -\frac{K}{\eta}$ for $i\in \NN$. From the definition of $\mathfrak{B}$ (see equation (\ref{Bell-a})), we have on one hand 
\begin{align*}
\inf_{a\in \mathbf{A}^{g}(x)} & \Big\{ \frac{K}{K+\eta}\int_{\mathbf{X}} V_{i}(\widetilde x) \widetilde{P}(d\widetilde x |x,a)+ \frac{1}{K+\eta} C^{g}(x,a) \Big\}
\geq -\frac{K}{K+\eta} \frac{K}{\eta} - \frac{K}{K+\eta}=-\frac{K}{\eta},
\end{align*}
and on the other hand
\begin{align*}
\inf_{a\in \mathbf{A}^{i}(x)}  \Big\{ \int_{\mathbf{X}} V_{i}(\widetilde x) Q(d\widetilde x | x,a) + c^{i}(x,a) \Big\} \geq \inf_{a\in \mathbf{A}^{i}(x)}  \Big\{ \int_{\mathbf{X}} V_{i}(\widetilde x) Q(d\widetilde x | x,a) \Big\} \geq 
-\frac{K}{\eta},
\end{align*}
since $c^{i}$ is non-negative (recalling Assumption \ref{HA4}).
Finally, combining the two previous equations, we obtain that $- \frac{K}{\eta} \leq V_{i+1}$.
Therefore, it follows that there exists a bounded function $V_{\infty}$ such that $V_{i}(x) \downarrow V_{\infty}(x)$ as $i\rightarrow \infty$, for any $x\in \mathbf{X}$ and so, $V_\infty\in \mathbb{U}(\mathbf{X})$ (see Theorem 4 in section 6, chapter 4 of \cite{bourbaki71}).
Now by using equation (\ref{ineqVi}), we obtain that $\mathfrak{B}V_{\infty}\leq \mathfrak{B}V_{n}\leq V_{n}$ for any $n\in\NN$ since the operator $\mathfrak{B}$ is monotone.
This implies that $\mathfrak{B}V_{\infty}\leq V_{\infty}$. Again from  (\ref{ineqVi}), it follows that
$V_{\infty} \leq \mathfrak{B} V_{i}$ for any $i\in \NN$.
Consequently, for $x\in \mathbb{X}^{i}$ and for any $a^g\in \mathbf{A}^{g}(x)$ and $a^i\in \mathbf{A}^{i}(x)$ we have
\begin{align*}
V_{\infty} (x) \leq  \Bigg\{ \frac{K}{K+\eta}\int_{\mathbf{X}} V_{i}(y) \widetilde{P}(dy |x,a^g)+ \frac{1}{K+\eta} C^{g}(x,a^g) \Bigg\} 
\wedge \Bigg\{ \int_{\mathbf{X}} V_{i}(y) Q(dy | x,a^i) + c^{i}(x,a^i) \Bigg\}.
\end{align*}
Now by taking the limit as $i\rightarrow \infty$ in the previous equation and by using the bounded convergence Theorem, it yields
\begin{align*}
V_{\infty} (x) \leq  \Bigg\{ \frac{K}{K+\eta}\int_{\mathbf{X}} V_{\infty}(y) \widetilde{P}(dy |x,a^g)+ \frac{1}{K+\eta} C^{g}(x,a^g) \Bigg\} 
\wedge \Bigg\{ \int_{\mathbf{X}} V_{\infty}(y) Q(dy | x,a^i) + c^{i}(x,a^i) \Bigg\},
\end{align*}
showing that $V_{\infty}(x) \leq \mathfrak{B} V_{\infty}(x)$. By using the similar arguments, it is easy to show that the case where $x\notin \mathbb{X}^{i}$ leads to the same conclusion, that is $V_{\infty}(x) \leq \mathfrak{B} V_{\infty}(x)$ for all $x\in\mathbf{X}$.
Finally, we have shown that $V_{\infty}=\mathfrak{B} V_{\infty}$ and so function $V$ defined by $V_\infty$ solves equation (\ref{Bell}).

The set $\mathbf{X}^g$ coincides with 
\begin{equation}\label{edot}
\left\{x\in\mathbf{X}: V(x)=\inf_{a\in \mathbf{A}^{g}(x)}  \Big\{ \frac{K}{K+\eta}\int_{\mathbf{X}} V(\widetilde x) \widetilde{P}(d\widetilde x |x,a)+ \frac{1}{K+\eta} C^{g}(x,a) \Big\}\right\}
\end{equation}
and hence is Borel-measurable, as well as $\mathbf{X}^i$.

According to Proposition 7.34 in \cite{bertsekas78}, for any $\varepsilon>0$, there is a measurable map $\widetilde\varphi^g$ from $\bf X$ to $\mathbf{A}^g$ such that for all $z\in\mathbf{X}$
\begin{align*}
\widetilde\varphi^g(z) \in \Bigg\{  \widetilde a\in\mathbf{A}^g(z): & \frac{K}{K+\eta}\int_{\mathbf{X}} V(\widetilde x) \widetilde{P}(d\widetilde x |x,\widetilde a)+ \frac{1}{K+\eta} C^{g}(x,\widetilde a) \\
& \leq \inf_{a\in \mathbf{A}^{g}(x)}  \Big\{ \frac{K}{K+\eta}\int_{\mathbf{X}} V(\widetilde x) \widetilde{P}(d\widetilde x |x,a)+ \frac{1}{K+\eta} C^{g}(x,a) \Big\}+\frac{\varepsilon}{K+\eta} \Bigg\}
\end{align*}
The restriction of  $\widetilde\varphi^g$ to $\mathbf{X}^g$ gives the mapping  $\varphi^g$ as required.
The mapping $\varphi^i$ is built in a similar way, working in the space $\mathbb{X}^{i}$ and passing to $\mathbf{X}^i$.
This gives the result.
\hfill $\Box$

\bigskip

\begin{proposition}
\label{existence-Bell2} 
Suppose Assumptions \ref{HA} and \ref{HC2} hold.
Then the increasing sequence of functions 
$(V_{i})_{i\in \NN}$ defined iteratively by $V_{i+1}=\mathfrak{B}V_{i}$ with $V_{0}=-\frac{K}{\eta}$
belongs to $\mathbb{B}(\mathbf{X})$ and converges to a bounded lower semicontinuous function $V$ on $\mathbf{X}$ satisfying the Bellman equation (\ref{Bell}). Moreover, the corresponding sets $\mathbf{X}^g$ and $\mathbf{X}^i$ are measurable
and, for any $\varepsilon\geq0$, there exists a Borel-measurable mapping $\varphi^i:\mathbf{X}^i\to\mathbf{A}^i$ 
(respectively, $\varphi^g:\mathbf{X}^g\to\mathbf{A}^g$) satisfying equation (\ref{qu1}) (respectively, equation (\ref{qu2})).
\end{proposition} 
\textbf{Proof:}
According to Proposition 7.33 in \cite{bertsekas78} and by considering the sequence $(V_{i})_{i\in \NN}$ defined by $V_{i+1}=\mathfrak{B} V_{i}$  with $V_{0}=-\frac{K}{\eta}$, it can be shown by using the same arguments as in Proposition \ref{existence-Bell} that $|V_i(x)|\le K/\eta$, $V_{i+1}\ge V_i$ and
$V_i$ is lower semicontinuous  for any $i\in \NN$. Consequently, $(V_{i})_{i\in \NN}$ converges pointwise to a limit denoted by  $V_\infty$ which is lower semicontinuous.

Clearly, $\mathfrak{B} V_\infty\ge\mathfrak{B} V_n\ge V_n$ for any $n\in\NN$, so that $\mathfrak{B} V_\infty\ge V_\infty$.
To show the reverse inequality, consider a sequence $(\wb{a}_{i})_{i\in \NN}$ of measurable mappings from $\mathbf{X}$ to $\mathbf{A}$ satisfying
$\wb{a}_{i}(x)\in \mathbf{A}^{g}(x)\union \mathbf{A}^{i}(x)$ and reaching the infimum in 
\begin{align}
\mathfrak{B}V_{i}(x)=\inf_{a\in \mathbf{A}^{g}(x)} & \Big\{ \frac{K}{K+\eta}\int_{\mathbf{X}} V_{i}(\widetilde x) \widetilde{P}(d\widetilde x |x,a)+ \frac{1}{K+\eta} C^{g}(x,a) \Big\} \nonumber \\
& \wedge \inf_{a\in \mathbf{A}^{i}(x)}  \Big\{ \int_{\mathbf{X}} V_{i}(\widetilde x) Q(d\widetilde x | x,a) + c^{i}(x,a) \Big\},
\label{aiBell}
\end{align}
for any $x\in\mathbf{X}$.
Fix an arbitrary $x\in\mathbf{X}$. There exists a subsequence $(\wb{a}_{i_{j}}(x))_{j\in \NN}$ of $(\wb{a}_{i}(x))_{i\in \NN}$ that belongs either to $\mathbf{A}^g(x)$ or $\mathbf{A}^i(x)$.
Consider that $\wb{a}_{i_{j}}(x)\in \mathbf{A}^g(x)$ for any $j\in \NN$ (the other possibility can be dealt by using the same arguments).
Moreover, there is no loss of generality to assume that this subsequence converges to some $\wb{a}\in \mathbf{A}^g(x)$ since $\mathbf{A}^g(x)$ is compact.
For $n\in \NN$ and $j \in\NN$ such that $n\leq i_{j}$, we have
\begin{align}
\frac{K}{K+\eta} & \int_{\mathbf{X}}  V_n(\widetilde x) \widetilde{P}(d\widetilde x |x,\widebar a_{i_j})+ \frac{1}{K+\eta} C^{g}(x,\widebar a_{i_j})
\nonumber \\
& \le  \frac{K}{K+\eta}\int_{\mathbf{X}} V_{i_j}(\widetilde x) \widetilde{P}(d\widetilde x |x,\widebar a_{i_j})+ \frac{1}{K+\eta} C^{g}(x,\widebar a_{i_j})
=\mathfrak{B}V_{i_j}(x)=V_{i_j+1}(x)\le V_\infty(x),
\label{alex1}
\end{align}
where the first inequality comes from the fact that $V_{i_j}\ge V_n$.
The real-valued mapping on $\mathbf{A}^g(x)$ defined by
$$\frac{K}{K+\eta}\int_{\mathbf{X}} V_n(\widetilde x) \widetilde{P}(d\widetilde x |x,.)+ \frac{1}{K+\eta} C^{g}(x,.)$$
is lower semicontinuous.
Consequently taking the limit as $j$ tends to infinity in equation (\ref{alex1}), it yields
$$\frac{K}{K+\eta}\int_{\mathbf{X}} V_n(\widetilde x) \widetilde{P}(d\widetilde x |x,\widebar a)+ \frac{1}{K+\eta} C^{g}(x,\widebar a)\le V_\infty(x).$$
The bounded convergence theorem implies
$$\frac{K}{K+\eta}\int_{\mathbf{X}} V_\infty(\widetilde x) \widetilde{P}(d\widetilde x |x,\widebar a)+ \frac{1}{K+\eta} C^{g}(x,\widebar a)\le V_\infty(x),$$
so that $\mathfrak{B}V_\infty(x)\le V_\infty(x)$, and hence $\mathfrak{B}V_\infty=V_\infty$ showing that the lower semicontinuous bounded function $V_\infty$ solves equation (\ref{Bell}).
The rest of the proof is similar to that of Proposition \ref{existence-Bell} but is now based on Proposition 7.33 in \cite{bertsekas78} and on the fact that
$V$ is lower semicontinuous.
\hfill $\Box$

\bigskip

\begin{remark}
\begin{itemize}
\item[(a)] The kernel $\frac{K}{K+\eta}\widetilde P$ is substochastic, so that the Bellman operator $\mathfrak{B}$ corresponds to a transient MDP. In this context, the convergence of the value iteration and the existence of the strict optimizer, (as $(\varphi^i,\varphi^g)$ for $\varepsilon=0$ according to our notation) has been established in \cite{james06} under conditions similar to \ref{HC2} but more restrictive.
\item[(b)] There are many other conditions guaranteeing similar statements of Propositions \ref{existence-Bell} of
\ref{existence-Bell2} to be correct. For example, instead of Condition \ref{HC1} one can require that
\begin{itemize}
\item the sets $\mathbf{A}^g$ and  $\mathbf{A}^i$  are $\sigma$-compact;
\item for any $F\in \mathbb{B}(\mathbf{X})$, the functions $\displaystyle
\int_\mathbf{X}F(z)\widebar q(dz|x,\cdot)$ and $\displaystyle \int_\mathbf{X}F(z)Q(dz|x,\cdot)$ are continuous for any $x\in\mathbf{X}$ and $x\in \mathbb{X}^{i}$ correspondingly;
\item the functions $C^g(x,\cdot)$ and $c^i(x,\cdot)$ are lower semicontinuous for any $x\in\mathbf{X}$ and $x\in \mathbb{X}^{i}$ correspondingly.
\end{itemize}
This condition will ensure
the existence of a function $V\in \mathbb{B}(\mathbf{X})$ satisfying the Bellman equation (\ref{Bell}), the measurability of the corresponding sets $\mathbf{X}^g$ and $\mathbf{X}^i$, and the existence of a Borel-measurable mapping $\varphi^i:\mathbf{X}^i\to\mathbf{A}^i$ 
(respectively, $\varphi^g:\mathbf{X}^g\to\mathbf{A}^g$) satisfying equation (\ref{qu1}) (respectively, equation (\ref{qu2})).

If additionally $\mathbf{A}^g$ and $\mathbf{A}^i$ are compact then the formulae (\ref{qu1}) and (\ref{qu2}) hold also for $\varepsilon=0$
according to Corollary 4.3 in \cite{rieder78}.
\end{itemize}
\end{remark}

\bigskip

The next two technical lemmas are needed to construct optimal and $\varepsilon$-optimal control strategies in Theorem \ref{th1}.
\begin{lemma}
\label{Ineq-1}
Suppose Assumption \ref{HA} and either of the assumptions \ref{HC1} or \ref{HC2} hold. Let $V$ be a bounded measurable solution of the Bellman equation (\ref{Bell}). Then the following assertions hold.
\begin{itemize}
\item[(a)] For any $x\in \mathbf{X}$ and $b\in \Xi$,
\begin{eqnarray*}
\int_{\mathbf{Y}} \Big[ C^{i}(y) + V(\widebar{x}(y) \Big] \beta^{b}(dy|x) & \geq & V(x).
\end{eqnarray*}
\item[(b)] If additionally Assumption \ref{HB} holds then, for any $\varepsilon>0$, there is a Markov non-randomized policy $b^*\in\Xi$ for the controlled model
$\mathcal{M}^{i}$ such that for any $x\in \mathbf{X}$
\begin{eqnarray}
\label{cb*1}
\int_{\mathbf{Y}} \Big[ C^{i}(y) + V(\widebar{x}(y) \Big] \beta^{b^*}(dy|x)  \leq  V(x)+\varepsilon
\end{eqnarray}
and
\begin{eqnarray}
\label{cb*2}
\beta^{b^*}\big( \{y\in Y:\widebar{x}(y)\in \mathbf{X}^g\} |x \big)=1.
\end{eqnarray}
\end{itemize}
Moreover, under Assumptions \ref{HA}, \ref{HB} and \ref{HC2}, the statement of item (b) can be strengthened. Indeed, it holds for $\varepsilon=0$ for a stationary non-randomized policy $b^*\in\Xi$.
\end{lemma}
\textbf{Proof:}
Associated with the discrete-time MDP $\mathcal{M}^{i}$,
consider the cost per stage function defined on $\mathbf{X}_{\Delta}\times \mathbf{A}^{i}_{\Delta}$ by $D=c^{i}+I_{\mathbf{X}\times \{\Delta\}} V$.
Let $x\in \mathbf{X}$ and $b$ be an arbitrary policy for $\mathcal{M}^{i}$
generating the process $(\widetilde{x}_{j},\wt{a}_{j})_{j\in \NN}$
with initial distribution $\delta_{x}$ and the corresponding strategic measure
$\beta^{b}(\cdot|x)$  on $\big( \mathbf{X}_{\Delta} \times \mathbf{A}^{i}_{\Delta} \big)^{\infty}$.
$E_{x}^{b}[\cdot]$ represents the expectation with respect to this strategic measure $\beta^{b}(\cdot|x)$.
Defining $\tau=\inf \{j\in \NN: \wt{a}_{j}=\Delta\}$, we obtain, by using the bounded convergence theorem and the definition of $Q_{\Delta}$,
\begin{eqnarray}
\lim_{m\rightarrow \infty} E_{x}^{b} \Big[\sum_{j=0}^{m} I_{\mathbf{X}\times \{\Delta\}}(\widetilde{x}_{j},\wt{a}_{j}) V(\widetilde{x}_{j}) \Big]
& = & E_{x}^{b} \Big[ \sum_{j=0}^{\infty} I_{\mathbf{X}\times \{\Delta\}}(\widetilde{x}_{j},\wt{a}_{j}) V(\widetilde{x}_{j}) \Big] \nonumber \\
& = & E^b_x \big[ V(\widetilde{x}_{\tau}) I_{\{\tau<\infty\}}\big].
\label{eq1}
\end{eqnarray}
Moreover, since $c^{i}$ is nonnegative
\begin{eqnarray}
\lim_{m\rightarrow \infty} E_{x}^{b} \Big[\sum_{j=0}^{m} c^{i}(\widetilde{x}_{j},\wt{a}_{j}) \Big]
= E_{x}^{b} \Big[\sum_{j=0}^{\infty} c^{i}(\widetilde{x}_{j},\wt{a}_{j}) \Big],
\label{eq2}
\end{eqnarray}
by the monotone convergence theorem.
Therefore, equations (\ref{eq1}) and (\ref{eq2}) yield
\begin{eqnarray}
\lim_{m\rightarrow \infty} E^{b}_{x}\Big[\sum_{j=0}^{m} D(\widetilde{x}_{j},\wt{a}_{j})\Big] & = & E^{b}_{x}\Big[\sum_{j=0}^{\infty} D(\widetilde{x}_{j},\wt{a}_{j})\Big].
\label{eq3}
\end{eqnarray}

\noindent
Regarding item $(a)$, we have $\beta^{b}(\{\tau<\infty\}|x)=1$ since $b\in \Xi$, and so
\begin{eqnarray}
E^{b}_{x}\Big[\sum_{j=0}^{\infty} D(\widetilde{x}_{j},\wt{a}_{j})\Big]=\int_{\mathbf{Y}} \big[ C^{i}(y) + V(\widebar{x}(y) \big] \beta^{b}(dy|x).
\label{eq4}
\end{eqnarray}

Consider the function $V_{\Delta}$ defined on $\mathbf{X}_{\Delta}$ by $V_{\Delta}(z)=V(z)$ if $z\in \mathbf{X}$ and $V_{\Delta}(\Delta)=0$.
From equation (\ref{Bell}) it is easy to show that $V_{\Delta}$ satisfies the following inequality
\begin{eqnarray*}
 \inf_{a\in \mathbf{A}^{i}_{\Delta}(x)}  \Big\{ D(x,a) + \int_{\mathbf{X}} V_{\Delta}(z) Q_\Delta(dz | x,a) \Big\} & \geq & V_{\Delta}(x), \quad x\in \mathbf{X}_\Delta.
\end{eqnarray*}
We have
$$E^{b}_{x}[V_{\Delta}(\widetilde{x}_{m+1})|\sigma\{(\widetilde{x}_{j},\wt{a}_{j}): j\in \NN_{m}\}]=\int_{\mathbf{X}_{\Delta}} V_{\Delta}(z)
Q_{\Delta}(dz|\widetilde{x}_{m},\wt{a}_{m})\geq V_{\Delta}(\widetilde{x}_{m})-D(\widetilde{x}_{m},\wt{a}_{m}),$$
for any $m\in \NN$.
Since $V_{\Delta}$ is bounded, it yields that
$$E^{b}_{x}\Big[\sum_{j=0}^{m} D(\widetilde{x}_{j},\wt{a}_{j})\Big]\geq V_{\Delta}(x)-E^{b}_{x}[V_{\Delta}(\widetilde{x}_{m+1})]
=V(x)-E^{b}_{x}[V_{\Delta}(\widetilde{x}_{m+1})],$$
for any $m\in \NN$.
Therefore, taking the limit as $m$ tends to infinity in the previous inequality and by using equations (\ref{eq3})-(\ref{eq4}) we obtain
$$\int_{\mathbf{Y}} \big[ C^{i}(y) + V(\widebar{x}(y) \big] \beta^{b}(dy|x) \geq V(x)
- \limsup_{m\rightarrow \infty}E^{b}_{x} \big[V_{\Delta}(\widetilde{x}_{m+1})\big].$$
However, $\big| E^{b}_{x}[V_{\Delta}(\widetilde{x}_{m})] \big| \leq \sup_{z\in X}|V(z)| \beta^{b}(\{m\leq \tau\}|x)$ and so,
$\limsup_{m\rightarrow \infty}E^{b}_{x}[V_{\Delta}(\widetilde{x}_{m+1})]=0$ since $\beta^{b}(\{\tau=\infty\}|x)=0$,  showing the result.

\bigskip

\noindent
To prove item $(b)$, introduce the following Markov non-randomized policy $b^*$ for the controlled model $\mathcal{M}^{i}$ defined by
$b^*=\big(\varphi^{i}_{j}\big)_{j\in\NN}$ where for $j\in \NN$, $\varphi^{i}_{j}$ is the $\mathbf{A}^i_{\Delta}$-valued measurable mapping defined on $\mathbf{X}_{\Delta}$ satisfying the following requirements:
\begin{itemize}
\item if $x\in \mathbf{X}^{i}$ then $\varphi^{i}_{j}(x)$ is an element of the set
$$\left\{a\in A^i(x) : \int_{\mathbf{X}} V(z)Q(dz|x,a)+c^i(x,a)\le V(x)+\varepsilon\left(\frac{1}{2}\right)^{j+1}\right\},$$
which is not empty by the definition of $ \mathbf{X}^{i}$.
\item if $x\in \mathbf{X}^{g}\union\{\Delta\}$ then $\varphi^{i}_{j}(x)=\Delta$.
\end{itemize}
The existence of such a measurable mapping was established in Proposition \ref{existence-Bell} under Assumptions
\ref{HA} and \ref{HC1} and in Proposition \ref{existence-Bell2} under Assumptions \ref{HA} and \ref{HC2}.
Now, for any $j\in \NN$ and $x\in \mathbf{X}_\Delta$,
$$D(x,\varphi^i_j(x))+\int_{\mathbf{X}_\Delta}V_\Delta(z)Q_{\Delta}(dz|x,\varphi^i_j(x))\le V_\Delta(x)+\varepsilon\left(\frac{1}{2}\right)^{j+1}.$$
Indeed, the previous inequality clearly holds for $x\in \mathbf{X}^{g}\union\{\Delta\}$. Now, if $x\in \mathbf{X}^{i}$ then it follows from the definition of
$\varphi^i_j(x)$.
Consequently,  for any $m\in \NN$ we have
\begin{eqnarray*}
E^{b^{*}}_{x}[V_{\Delta}(\widetilde{x}_{m+1})|\sigma\{(\widetilde{x}_{j},\wt{a}_{j}): j\in \NN_{m}\}]&=&\int_{\mathbf{X}_{\Delta}} V_{\Delta}(z)
Q(dz|\widetilde{x}_{m},\wt{a}_{m})\\
&\leq& V_{\Delta}(\widetilde{x}_{m})-D(\widetilde{x}_{m},\wt{a}_{m})+\varepsilon\left(\frac{1}{2}\right)^{j+1}.
\end{eqnarray*}
Therefore, by using the fact that $V_{\Delta}$ is bounded, we obtain 
$$E^{b^{*}}_{x}\Big[\sum_{j=0}^{m} D(\widetilde{x}_{j},\wt{a}_{j})\Big]
\leq V_{\Delta}(x)-E^{b^*}_{x}[V_{\Delta}(\widetilde{x}_{m+1})]+
\varepsilon\cdot\frac{1}{2}\cdot\frac{1-\left(\frac{1}{2}\right)^{m+1}}{1-\frac{1}{2}},$$
implying that
$\limsup_{m\to\infty}E^{b^*}_{x}\Big[\sum_{j=0}^{m} 
c^{i}(\widetilde{x}_{j},\wt{a}_{j})\Big]$ is finite.
Moreover, observe that
$$\{\tau=\infty\} \subset\Big\{ \limsup_{m\to\infty} \sum_{j=0}^{m} c^{i}(\widetilde{x}_{j},\wt{a}_{j}) =\infty \Big\}$$
by using Assumption \ref{HB}.
By the monotone convergence theorem, it yields $\beta^{b^{*}}\big( \{ \tau <\infty \}|x\big)=1$ implying $b^*\in\Xi$.
Now, by using similar arguments as for the proof of item $(a)$, we get
$$\int_{Y} \Big[ C^{i}(y) + V(\widebar{x}(y) \Big] \beta^{b^*}(dy|x) = \lim_{m\to\infty} E^{b^*}_x \Big[\sum_{j=0}^m D(\widetilde x_j,\wt{a}_{j})\Big]
\leq  V_\Delta(x)+\varepsilon=V(x)+\varepsilon.$$
Finally, observe that $\{ \tau <\infty \} \subset \big\{y\in \mathbf{Y}:\widebar{x}(y)\in \mathbf{X}^g\}$ and so we get the last assertion.

\bigskip

The proof of the last statement of Lemma \ref{Ineq-1} is similar to part (b), with a reference to Proposition \ref{existence-Bell2}.
\hfill $\Box$

\bigskip

\begin{lemma}
\label{Ineq-3}
Suppose Assumption \ref{HA} and either of the assumptions \ref{HC1} or \ref{HC2} hold.
Let $x\in \mathbf{X}$, $\beta\in \mathcal{P}^{\mathbf{Y}}(x)$ and $V$ be a bounded solution of the Bellman equation (\ref{Bell}). Then
\begin{align*}
-\eta V(x)  +\int_{\mathbf{X}} \int_{\mathbf{Y}} \big[V(\widebar{x}(y)) +C^{i}(y) \big] \beta(dy | z) \widebar{q}(dz| x,a)
-V(x) \widebar{q}(\mathbf{X}| x,a) +C^{g}(x,a) \geq 0,
\end{align*}
for any $(x,a)\in \mathbb{K}^{g}$.
\end{lemma}
\textbf{Proof:}
From Lemma \ref{Ineq-1}, we have that for any $z\in \mathbf{X}$
\begin{eqnarray}
\int_{\mathbf{Y}} \big[V(\widebar{x}(y)) +C^{i}(y) \big] \beta(dy | z) \geq V(z),
\end{eqnarray}
and so recalling that $\widebar{q}$ is a positive kernel
\begin{align*}
-\eta V(x)  & +\int_{\mathbf{X}} \int_{\mathbf{Y}} \big[V(\widebar{x}(y)) +C^{i}(y) \big] \beta(dy | z) \widebar{q}(dz| x,a)
-V(x) \widebar{q}(\mathbf{X}| x,a) +C^{g}(x,a) \nonumber \\
& \geq -\eta V(x)  +\int_{\mathbf{X}} V(z) \widebar{q}(dz| x,a) -V(x) \widebar{q}(\mathbf{X}| x,a) +C^{g}(x,a)
\end{align*}
for any $(x,a)\in \mathbb{K}^{g}$. Now the result follows from equation (\ref{Bell}).
\hfill $\Box$

\bigskip

The next result shows the existence of optimal and $\varepsilon$-optimal control strategies.
\begin{theorem}\label{th1}
Suppose assumptions \ref{HA}, \ref{HB} hold. Let $V$ be a bounded measurable  solution of the Bellman equation (\ref{Bell}). 
\begin{itemize}
\item[(a)] If assumptions \ref{HC1} or \ref{HC2} hold, then for any control strategy 
$u=(u_{n})_{n\in\NN}\in \mathcal{U}$ with
$u_{n}=\big( \psi_{n},\pi_{n},\gamma^0_{n},\gamma^1_{n} \big)$ for $n\in \NN^{*}$
$$\mathcal{V}(u,x_{0})\ge V(x_0)$$
and, for any $\varepsilon>0$, there is a non-randomized almost stationary strategy $u^*$ such that
$$\mathcal{V}(u^*,x_{0})\le V(x_0)+\varepsilon,$$
which satisfies $\psi_n(\cdot|h_n)=\delta_{\infty}(\cdot)$, that is, the interventions occur only after the natural jumps (and maybe at the initial time moment $t=0$). 
\item[(b)] If assumption \ref{HC2} holds then there exists a non-randomized stationary strategy $u^*$ such that 
$$\mathcal{V}(u^*,x_{0})= V(x_0),$$ which satisfies $\psi_n(\cdot|h_n)=\delta_{\infty}(\cdot)$, that is, the interventions occur only after the natural jumps (and maybe at the initial time moment $t=0$).
\end{itemize}
\end{theorem}
\textbf{Proof:} 
(a) If $\mathcal{V}(u,x_{0})=+\infty$ then the inequality is clearly satisfied.
Now, if $\mathcal{V}(u,x_{0})<+\infty$, then  we have $\mathbb{P}^{u}_{x_{0}}(T_{\infty}<\infty)=0$ from Proposition \ref{explo}.
Consequently, Lemma \ref{lem1} applied to a bounded solution $V$ of the Bellman equation (\ref{Bell}) yields
\begin{align}
& \mathcal{V}(u,x_{0}) =  \int_{\mathbf{Y}} V(\widebar{x}(y)) u_{0}(dy|x_{0}) +\int_{\mathbf{Y}} C^{i}(y)  u_{0}(dy|x_{0})  \nonumber \\
&+ \mathbb{E}^{u}_{x_{0}} \Bigg[ \int_{0}^{+\infty} \int_{\mathbf{A}^{g}}
e^{-\eta s} \Big[ -\eta V(\wb{x}(\xi_{s}))+C^{g}(\wb{x}(\xi_{s-}),a) \nonumber \\
& +\int_{\mathbf{X}} \int_{\mathbf{Y}} \big\{ V(\widebar{x}(y)) + C^{i}(y) \big\} \gamma^{0}(dy | v,s) \widebar{q}(dv| \wb{x}(\xi_{s}),a)
-V(\wb{x}(\xi_{s})) \widebar{q}(\mathbf{X}| \wb{x}(\xi_{s}),a)    \Big] \pi(da |s) ds \Bigg] \nonumber \\
& + \mathbb{E}^{u}_{x_{0}} \Bigg[ \sum_{n\in \NN^{*}}\int_{]T_{n}, T_{n+1}] } \int_{\mathbf{Y}} e^{-\eta s}
\Big[ V(\widebar{x}(y))+C^{i}(y) - V(\wb{x}(\xi_{s-})) \Big]  
\gamma^{1}_{n}(dy |H_{n}) \frac{\psi_{n}(ds-T_{n} | H_{n})}{\psi_{n}([s-T_{n},+\infty] | H_{n})} \Bigg].\label{eq18}
\end{align}
Observe that, $\gamma^{1}_{n}(dy |h_{n}) \in \mathcal{P}^{\mathbf{Y}^{*}}({\wb{x}(y_n)})
\subset\mathcal{P}^{\mathbf{Y}}({\wb{x}(y_n)})$ and
$\gamma^{0}_{n}(dy |h_{n}, s,x) \in \mathcal{P}^{\mathbf{Y}}(x)$ for any $n\in\NN^{*}$, $s\in \RR_{+}$, $x\in \mathbf{X}$
and $h_{n}=(y_0,\theta_1,y_1,\ldots,\theta_{n},y_{n})\in \mathbf{H}_{n}$.
Consequently, we obtain from Lemma \ref{Ineq-1} that
\begin{eqnarray*}
\int_{\mathbf{Y}} \Big[ V(\widebar{x}(y)) +C^{i}(y) - V(\wb{x}(y_{n}))\Big] \gamma^{1}_{n}(dy |h_{n}) & \geq & 0
\end{eqnarray*}
for any $n\in\NN^{*}$ and $h_{n}=(y_0,\theta_1,y_1,\ldots,\theta_{n},y_{n})\in \mathbf{H}_{n}$.
Now, by recalling Lemma \ref{Ineq-3}, we have that
\begin{align*}
-\eta V(\wb{x}(y_{n}))  +\int_{\mathbf{X}} \int_{\mathbf{Y}} \big[V(\widebar{x}(y)) +C^{i} & (y) \big] \gamma^{0}_{n}(dy |h_{n}, s,x)
\widebar{q}(dx| \wb{x}(y_{n}),a)
\nonumber \\
& -V(\wb{x}(y_{n})) \widebar{q}(\mathbf{X}| \wb{x}(y_{n}),a) +C^{g}(\wb{x}(y_{n}),a) \geq 0
\end{align*}
for any $n\in\NN^{*}$, $s\in \RR_{+}$, $h_{n}=(y_0,\theta_1,y_1,\ldots,\theta_{n},y_{n})\in \mathbf{H}_{n}$ and $a\in \mathbf{A}^{g}(\wb{x}(y_{n}))$.
Observe that $\xi_{s-}=Y_{n}$ on the stochastic interval  $\sr T_{n}, T_{n+1}\sr$ and that
$\xi_{s-}=\xi_{s}$ on stochastic interval  $\sr T_{n}, T_{n+1} \sl$.
Therefore, the two previous equations yields
\begin{align}
& \int_{0}^{+\infty} \int_{\mathbf{A}^{g}}
e^{-\eta s} \Big[ -\eta V(\wb{x}(\xi_{s}))+C^{g}(\wb{x}(\xi_{s-}),a)
+\int_{\mathbf{X}} \int_{\mathbf{Y}} \big\{ V(\widebar{x}(y)) + C^{i}(y) \big\} \gamma^{0}(dy |x,s) \widebar{q}(dx| \wb{x}(\xi_{s}),a)
\nonumber \\
& -V(\wb{x}(\xi_{s})) \widebar{q}(\mathbf{X}| \wb{x}(\xi_{s}),a)    \Big] \pi(da |s) ds
+ \sum_{n\in \NN^{*}}\int_{]T_{n}, T_{n+1}] } \int_{ \mathbf{Y}} e^{-\eta s}
\Big[ V(\widebar{x}(y))+C^{i}(y) - V(\wb{x}(\xi_{s-})) \Big]  \nonumber \\
& \gamma^{1}_{n}(dy |H_{n}) \frac{\psi_{n}(ds-T_{n} | H_{n})}{\psi_{n}([s-T_{n},+\infty] | H_{n})} \geq 0,
\quad \text{ $\mathbb{P}^{u}_{x_{0}}-$ a.s.}
\end{align}
implying that
\begin{eqnarray*}
\mathcal{V}(u,x_{0}) &  \geq &  \int_{\mathbf{Y}} V(\widebar{x}(y)) u_{0}(dy|x_{0}) +\int_{\mathbf{Y}} C^{i}(y)  u_{0}(dy|x_{0}).
\end{eqnarray*}
Finally, we get the first assertion, that is  $\mathcal{V}(u,x_{0}) \geq V(x_0)$, by using item (a) of Lemma \ref{Ineq-1} and recalling that $u_0\in{\cal P}^\mathbf{Y}(x_0)$. 

\bigskip

Fix an arbitrary $\varepsilon>0$ and consider the following non-randomized almost stationary control strategy $u^*=(u^*_n)_{n\in\NN}$
with $u^{*}_{n}=\big( \psi_{n},\pi_{n},\gamma^0_{n},\gamma^1_{n} \big)$ for $n\in \NN^{*}$, given by the following elements.
Let $$\psi_n(\cdot|h_n)=\delta_\infty(\cdot)$$ (i.e. the interventions occur only after the natural jumps and maybe at the initial time moment $t=0$).
Set $$\pi_n(\cdot|h_n,t)=\delta_{\varphi^g(\wb{x}(y_n))}(\cdot),$$ where, for $\widetilde{x}\in \mathbf{X}^g$,
$$\displaystyle\varphi^g(\widetilde x)\in\left\{a\in \mathbf{A}^g(\widetilde x):~\int_\mathbf{X} V(v)\widebar q(dv|\widetilde x,a)-V(\widetilde x)
\widebar q(\mathbf{X}|\widetilde x,a)+C^g(\widetilde x,a)\le \eta V(\widetilde x)+\frac{\eta\varepsilon}{3}\right\} $$
and, for $\widetilde x\in \mathbf{X}^i$, $\varphi^g(\cdot)$ is an arbitrary measurable mapping from $\mathbf{X}^i$ to $\mathbf{A}^g$ with $\varphi^g(x)\in\mathbf{A}^g(x)$. The existence of such a mapping follows from Proposition \ref{existence-Bell} under Assumptions \ref{HA} and \ref{HC1} and from Proposition \ref{existence-Bell2} under Assumptions \ref{HA} and \ref{HC2}.
Let $$\gamma^0_n(\cdot|h_n,\widetilde x)=\beta^{b^*}(\cdot|\widetilde x),$$
with the policy $b^*\in\Xi$ introduced in item (b) of  Lemma \ref{Ineq-1} and satisfying the following inequality
\begin{eqnarray}
\int_{\bf Y}[C^i(y)+V(\widebar x(y))]\beta^{b^*}(dy| x)\le V(x)+\min\left\{1,\frac{\eta}{K}\right\}\frac{\varepsilon}{3}
\label{betabcont1}
\end{eqnarray}
and  equation (\ref{cb*2}) for any $x\in X$.\\
Consider $\gamma^1_{n}$ as an arbitrary stochastic kernel on $\mathbf{Y}$  given $ \mathbf{H}_{n}$
satisfying $\gamma^1_{n}(\cdot|h_{n})\in  \mathcal{P}^{\mathbf{Y}^{*}}(\wb{x}(y_{n}))$
for any $h_{n}=(y_0,\theta_1,\ldots \theta_{n},y_{n})\in \mathbf{H}_{n}$
with $\wb{x}(y_{n})\in \mathbb{X}^{i}$.
Finally, set $u_0=\beta^{b^*}(\cdot|x_0)$.

\bigskip

\noindent
First, from Remark \ref{Cond-dist} and the definition of $G_{n}$ (see equation (\ref{defG2})) it follows that
\begin{eqnarray*}
\mathbb{P}^{u^{*}}_{x_{0}} \big(\Theta_{n+1} \in \Gamma_{\theta} | \mathcal{F}_{T_{n}}\big) & = & \int_{\Gamma_{\theta}} \lambda_{n}(\mathbf{X},H_{n},t)
e^{-\Lambda_{n}(\mathbf{X},H_{n},t)} dt.
\end{eqnarray*}
However, Assumption \ref{HA1} ensures that $\lambda_{n}(\mathbf{X},H_{n},t)$ is uniformly bounded by constant $K$ and so,
$\mathbb{P}^{u^{*}}_{x_{0}}(T_{\infty}<\infty)=0$.
Consequently, according to (\ref{eq18}) and the definition of the strategy $u^*$, we have
\begin{align*}
 \mathcal{V}(u^{*},x_{0}) =  \int_{\mathbf{Y}} V(\widebar{x}(y)) \beta^{b^*}(dy|x_{0}) & +\int_{\mathbf{Y}} C^{i}(y)  \beta^{b^*}(dy|x_{0})  \nonumber \\
+ \mathbb{E}^{u^{*}}_{x_{0}} \Bigg[ \int_{0}^{+\infty}  e^{-\eta s} \Big[ & -\eta V(\wb{x}(\xi_{s}))+C^{g}(\wb{x}(\xi_{s-}),\varphi^{g}(\wb{x}(\xi_{s-}))) \Big] ds \Bigg]  \nonumber \\
+ \mathbb{E}^{u^{*}}_{x_{0}} \Bigg[ \int_{0}^{+\infty}  e^{-\eta s}
\Big[ & \int_{\mathbf{X}} \int_{\mathbf{Y}} \big\{ V(\widebar{x}(y)) + C^{i}(y) \big\} \beta^{b^*}(dy | x) \widebar{q}(dx| \wb{x}(\xi_{s}),\varphi^{g}(\wb{x}(\xi_{s-})))
\nonumber \\
& -V(\wb{x}(\xi_{s})) \widebar{q}(\mathbf{X}| \wb{x}(\xi_{s}),\varphi^{g}(\wb{x}(\xi_{s-})))    \Big] ds \Bigg].
\end{align*}
Now, from inequality (\ref{betabcont1}) we have
\begin{eqnarray}
 \mathcal{V}(u^{*},x_{0}) & \leq & V(x_0)+\frac{\varepsilon}{3} 
+   \frac{\eta \varepsilon}{3K} \mathbb{E}^{u}_{x_{0}} \Big[ \int_{0}^{+\infty}  e^{-\eta s}  \widebar{q}(\mathbf{X}| \wb{x}(\xi_{s}),\varphi^{g}(\wb{x}(\xi_{s-}))) ds \Big]
\nonumber \\
&& + \mathbb{E}^{u}_{x_{0}} \Bigg[ \int_{0}^{+\infty}  e^{-\eta s} \mathcal{W}(s) ds \Bigg],
\label{af1}
\end{eqnarray}
where
\begin{eqnarray}
\mathcal{W}(s) & = & -\eta V(\wb{x}(\xi_{s}))+C^{g}(\wb{x}(\xi_{s-}),\varphi^{g}(\wb{x}(\xi_{s}))) \nonumber \\
& & + \int_{\mathbf{X}}V(x) \widebar{q}(dx| \wb{x}(\xi_{s}),\varphi^{g}(\wb{x}(\xi_{s})))
-V(\wb{x}(\xi_{s})) \widebar{q}(\mathbf{X}| \wb{x}(\xi_{s}),\varphi^{g}(\wb{x}(\xi_{s}))).
\label{defcalV}
\end{eqnarray}
Observe that since $V$ is bounded and by Assumption \ref{HA}, we have that $|\mathcal{W}|$ is uniformly bounded by a constant $K_{\mathcal{W}}$.
By Fubini's Theorem,
\begin{eqnarray}
\mathbb{E}^{u}_{x_{0}} \bigg[ \int_{0}^{+\infty}  e^{-\eta s} \mathcal W(s) ds \bigg] 
& \leq & \sum_{n\in\NN^{*}}\mathbb{E}^{u}_{x_{0}} \bigg[  I_{\mathbf{X}^{g}}(\wb{x}(Y_{n})) \int_{[T_{n}, T_{n+1}[}  e^{-\eta s} \mathcal{W}(s) ds \bigg]
\nonumber \\
& & + K_{\mathcal{W}} \sum_{n\in\NN^{*}}\mathbb{P}^{u}_{x_{0}} \Big(  \wb{x}(Y_{n}) \in \mathbf{X}^{i} \Big).
\label{af2}
\end{eqnarray}
Now, observe that 
\begin{eqnarray*}
\mathbb{P}^{u}_{x_{0}} \Big(  \wb{x}(Y_{n+1}) \in \mathbf{X}^{i} \big| \mathcal{F}_{T_{n}} \Big)
& = & G_{n} \Big( \Gamma^{i}_{y}\times \RR_{+} | H_{n} \Big),
\end{eqnarray*}
where $\Gamma^{i}_{y}=\{y\in \mathbf{Y} : \wb{x}(y)\in \mathbf{X}^{i}\}$.
However, $\gamma^{0}_{n}(\Gamma^{i}_{y}|H_{n},t,x)=\beta^{b^{*}}(\Gamma^{i}_{y}|x)=0$ for any $x\in \mathbf{X}$ according to equation (\ref{cb*2}), and so, from the definition of $G_{n}$,
we have $\mathbb{P}^{u}_{x_{0}} \Big(  \wb{x}(Y_{n+1}) \in \mathbf{X}^{i} \big| \mathcal{F}_{T_{n}} \Big)=0$ implying that
$\mathbb{P}^{u}_{x_{0}} \Big(  \wb{x}(Y_{n+1}) \in \mathbf{X}^{i} \Big)=0$ for any $n\in \NN^{*}$.
Moreover, $\mathbb{P}^{u}_{x_{0}} \Big(  \wb{x}(Y_{1}) \in \mathbf{X}^{i} \Big)=\beta^{b^{*}}(\Gamma^{i}_{y}|x_{0})$ according to equation (\ref{cb*2}).
Consequently, from equation (\ref{af2}), it follows that
\begin{eqnarray*}
\mathbb{E}^{u}_{x_{0}} \bigg[ \int_{0}^{+\infty}  e^{-\eta s} \mathcal{W}(s) ds \bigg] 
& \leq & \sum_{n\in\NN^{*}}\mathbb{E}^{u}_{x_{0}} \bigg[  I_{\mathbf{X}^{g}}(\wb{x}(Y_{n})) \int_{[T_{n}, T_{n+1}[}  e^{-\eta s} \mathcal{W}(s) ds \bigg],
\end{eqnarray*}
and so, from equation (\ref{defcalV}) and the definition of $\varphi^{g}$ on $\mathbf{X}^{g}$, it follows that 
\begin{eqnarray}
\mathbb{E}^{u}_{x_{0}} \bigg[ \int_{0}^{+\infty}  e^{-\eta s} \mathcal{W}(s) ds \bigg] 
& \leq & \frac{\varepsilon}{3}.
\label{af3}
\end{eqnarray}
Finally, combining equations (\ref{af1}), (\ref{af3}) and Assumption \ref{HA1}, we get the result $$\mathcal{V}(u^*,x_0)\leq V(x_{0})+\varepsilon.$$

(b) The proof is similar to the proof of item (a) by using now Proposition \ref{existence-Bell2} and the last statement of Lemma \ref{Ineq-1}.
\hfill $\Box$

\bigskip
The next result shows the existence of uniformly optimal and $\varepsilon$-optimal control strategies.
\begin{corollary} 
\label{corol1}
The following assertions hold.
\begin{itemize}
\item[(a)]  If assumptions \ref{HA}, \ref{HB} and \ref{HC1} hold then, for any $\varepsilon>0$, there exists a non-randomized almost stationary uniformly $\varepsilon$-optimal strategy $u^*$  satisfying $\psi_n(\cdot|h_n)=\delta_{\infty}(\cdot)$, that is, the interventions occur only after the natural jumps (and maybe at the initial time moment $t=0$).
\item[(b)]  If assumptions \ref{HA}, \ref{HB} and \ref{HC2} hold then there is a non-randomized stationary uniformly optimal strategy $u^*$  satisfying $\psi_n(\cdot|h_n)=\delta_{\infty}(\cdot)$, that is, the interventions occur only after the natural jumps (and maybe at the initial time moment $t=0$).
\item[(c)] In either case, 
$$\inf_{u\in{\cal U}}\mathcal{V}(u,x_{0})=V(x_0),$$ 
and the Bellman equation (\ref{Bell}) has a unique bounded measurable solution.
\end{itemize}
\end{corollary}
\textbf{Proof:} It follows directly from the proof of Theorem \ref{th1}: the control strategies $u^*$ does not depend on the initial state $x_0$.\hfill $\Box$

\bigskip

\begin{remark}
\label{decomp}
Roughly speaking, one should apply the gradual action $\varphi^g(x)$ if the current state is $x\in\mathbf{X}^g$, and one should apply the impulsive action $\varphi^i(x)$ if $x\in\mathbf{X}^i$. The mappings $\varphi^g$ and $\varphi^i$ have been introduced in Proposition  \ref{existence-Bell} and  \ref{existence-Bell2}
according to the Assumptions under consideration.
As rigorously described in the proof of Theorem \ref{th1}, the optimal (or $\varepsilon$-optimal) strategy $u^*$ looks as follows:
\begin{itemize}
\item $\psi_n(\cdot|h_n)=\delta_\infty(\cdot)$, i.e. the interventions occur only after the natural jumps (and maybe at the initial time moment $t=0$);
\item $\pi_n(\cdot|h_n,t)=\delta_{\varphi^g(\widebar x(y_n))}(\cdot)$, where, for $x\in\mathbf{X}^i$, $\varphi^g(\cdot)$ is an arbitrary fixed measurable mapping (it is in fact never applied);
\item $\gamma^0_n(\cdot|h_n,t,x)=\beta^{b^*}(\cdot|x)$, where $b^*\in\Xi$ is the Markov non-randomized policy in the model ${\cal M}^i$ defined by the mapping $\varphi^i_j:~\mathbf{X}^i\to\mathbf{A}^i$ complemented with equation $\varphi^i_j(x)=\Delta$ for $x\in\mathbf{X}^g$ ($j$ plays the role of time in the ${\cal M}^i$ model);
\item $\gamma^1_n(\cdot|h_n)$ is arbitrarily fixed (it is never used).
\end{itemize}
\end{remark}

It is obvious that usually one cannot avoid the series of instantaneous impulses when constructing an optimal control strategy because ${\cal V}(u,x_0)>V(x_0)$ if $\mathbb{P}^{u}_{x_{0}}(\exists t_1<t_2:\forall s\in(t_1,t_2)~\widebar x(\xi_s)\in\mathbf{X}^i)>0$. On the other hand, if one is looking for an $\varepsilon$-optimal strategy then instantaneous impulses can be avoided. Fix an $\frac{\varepsilon}{2}$-optimal strategy $u^*$ (see Theorem \ref{th1}(a)) and modify it in the following way. If a series of instantaneous impulses should be used, then one applies them sequentially, after small intervals $\delta_0,\delta_1,\ldots$. During those intervals, an arbitrarily fixed admissible gradual control is used. In case a natural jump occurs on these intervals, on the remainder time horizon one applies an arbitrary admissible gradual control, without any interventions. Denote this strategy as $\hat u$. Then $\hat u$ is $\varepsilon$-optimal, if the total length $\sum_{n\in\NN}\delta_n$ is small enough, e.g. if
  $$\sum_{n\in\NN}\delta_n\le\frac{\varepsilon}{4[K+K^2/\eta]},$$
where $K$ is the constant from Assumption \ref{HA}.

\section{Example}
\label{sec-example} 
We consider the following mathematical model of the epidemic with carriers. A closed population consists initially of $S$ susceptibles who can become infected (or ill) because of the contacts with carriers. The number of carriers changes according to the birth-and-death process, independent of the main population, with the birth and death rates $\rho_b(c)$ and $\rho_d(c)$ correspondingly, where $c$ is the current number of the carriers. Every one susceptible becomes infected at the rate $\kappa_i(c)$. Every one infective recovers and becomes immunized after an exponentially distributed time with parameter $\kappa_r$; the immunized individuals are removed from the consideration. The processes describing the dynamics of the number of carriers, the transformation of a susceptible to an infective, and the recovery of any one infective are mutually independent and uncontrolled. Thus, the set ${\bf A}^g$ consists of a single element meaning "do nothing", and we omit this argument. The transition rates from the current state $(s,c,i)$, where $s$
(respectively, $c$ and $i$) denotes the number of susceptibles (respectively, carriers and infectives), are as follows:
\begin{itemize}
\item $(s,c,i)\to (s,c\pm 1,i)$ at the rates $\rho_b(c)$ and $\rho_d(c)$ correspondingly; we assume that $\rho_b(0)=\rho_d(0)=0$;
\item $(s,c,i)\to (s-1,c,i+1)$ at the rate $s\cdot\kappa_i(c)$, and we assume that $\kappa_i(0)=0$;
\item $(s,c,i)\to (s,c,i-1)$ at the rate $i\cdot\kappa_r$.
\end{itemize}

The state space is ${\bf X}=\NN_S\times\NN\times\NN_{S+I}$; $x_0=(s_0=S,c_0,i_0=I)$ is the initial state, where $I$ is the initial number of infectives. The intensity of jumps $q(\tilde x|x)$ comes from the expressions given above in the obvious way. The cost rate $C^g(s,c,i)=i$ is associated with the infectives. Below we assume that all functions $\rho_b(c),\rho_d(c)$ and $\kappa_i(c)$ are bounded, non-negative.

The impulsive action is unique and means immunization of one susceptible: $\mathbb{X}^{i}=\NN_S^*\times\NN\times\NN_{S+I};$ $Q((s-1,c,i)|(s,c,i))=1$ and $c^i(x)=\lambda>0$ is the cost of one immunization. As previously, we omit the argument $a$ from all the formulae because the impulsive action is unique for all $x\in\mathbb{X}^{i}$.

Now the Bellman equation (\ref{Bell}) takes the form:
\begin{eqnarray}
&\min & \{ -\eta V(s,c,i)+\rho_b(c)V(s,c+1,i)+\rho_d(c)V(s,c-1,i)+s\kappa_i(c) V(s-1,c,i+1) \nonumber \\
&& +i\kappa_r V(s,c,i-1)-V(s,c,i)[\rho_b(c)+\rho_d(c)+s\kappa_i(c)+i\kappa_r]+i; \label{e32} \\
&& -V(s,c,i)+V(s-1,c,i)+\lambda\}=0, \nonumber
\end{eqnarray}
for $s>0$ and 
\begin{align*}
-\eta V(0,c,i)& +\rho_b(c)V(0,c+1,i)+\rho_d(c)V(0,c-1,i) \nonumber \\
& +i\kappa_r V(0,c,i-1)-V(0,c,i)[\rho_b(c)+\rho_d(c)+s\kappa_i(c)+i\kappa_r]+i=0,
\end{align*}
otherwise.
Here and below, we omit the brackets when indicating the state $x=(s,c,i)$.

All the assumptions \ref{HA}, \ref{HB} and \ref{HC} are satisfied if we fix the discrete topology in ${\bf X}$. Therefore, according to the Corollary \ref{corol1}, equation (\ref{e32}) has a single bounded solution.

If $s=0$ then $V(0,c,i)=\frac{i}{\eta+\kappa_r}$. Indeed, this function is bounded on $\NN\times \NN_{S+I}$ and solves equation
  $$V(0,c,i)=\frac{\rho_b(c) V(0,c+1,i)+\rho_d(c) V(0,c-1,i)+i\kappa_r V(0,c,i-1)+i}{\eta+\rho_b(c)+\rho_d(c)+i\kappa_r}.$$

\begin{lemma}\label{exl1} If $v$ is a bounded solution to equation
\begin{equation}\label{e33}
v(c)=\min\left\{\frac{\rho_b(c)v(c+1)+\rho_d(c)v(c-1)+\frac{\kappa_i(c)}{\eta+\kappa_r}}{\eta+\rho_b(c)+\rho_d(c)+\kappa_i(c)};~~\lambda\right\}~~c\in\NN,
\end{equation}
then $V(s,c,i)=sv(c)+\frac{i}{\eta+\kappa_r}$ is a bounded solution to (\ref{e32}) on $\bf X$. Morover, the first (respectively, second) expression in (\ref{e33}) equals $v(c)$ if and only if the first (respectively, second) expression in (\ref{e32}) is zero.
\end{lemma}
\textbf{Proof:} If we substitute expression $V(s,c,i)=sv(c)+\frac{i}{\eta+\kappa_r}$ into (\ref{e32}), we obtain the following equation:
$$\min\left\{ s\left[\rho_b(c)v(c+1)+\rho_d(c)v(c-1)+\frac{\kappa_i(c)}{\eta+\kappa_r}-v(c)[\eta+\rho_b(c)+\rho_d(c)+\kappa_i(c)]\right];~-v(c)+\lambda\right\}=0.$$

If $c$ is such that $v(c)=\lambda$ then the second expression in (\ref{e32}) equals zero and the first expression in (\ref{e32}) is non-negative
due to equation (\ref{e33}).
If $c$ is such that $v(c)$ equals the first expression in (\ref{e33}) then the first expression in (\ref{e32}) is zero and the second expression in (\ref{e32}) is non-negative due to equation (\ref{e33}). \hfill $\Box$

\bigskip

Let us derive the following technical result.
\begin{lemma}\label{exl2} Consider the non-negative, bounded, and increasing functions $\alpha_1$, $\alpha_2$ and $\alpha_3$ defined on $\NN$.
Suppose also that $\alpha_1(0)=\alpha_2(0)=\alpha_3(0)=0$ and  $\lim_{c\to\infty}(\alpha_1(c)+\alpha_2(c))<1$.
Consider $\lambda>0$ as a fixed constant. Then the following assertions hold.

(a) The equation
\begin{equation}\label{e34}
w(c)=\min\{\alpha_1(c)w(c+1)+\alpha_2(c)w(c-1)+\alpha_3(c);~~~~\lambda\}
\end{equation}
has a single bounded solution.

(b) Let $\lambda^*=\frac{\zeta_3}{1-\zeta_1-\zeta_2}$, where $\zeta_j=\lim_{c\to\infty}\alpha_j(c)$, $j=1,2,3$. \\
If $\lambda<\lambda^*$, then there is $c^*\in\NN^*$ such that the first term in the minimum in  (\ref{e34}) is strictly bigger than $\lambda$ if and only if $c\ge c^*$.\\
If $\lambda\ge\lambda^*$, then the first term in the minimum in (\ref{e34})  is not bigger than $\lambda$ and hence equals $w(c)$ for all $c\in\NN$. (In this case we say that $c^*=+\infty$.)
\end{lemma}

\textbf{Proof:} (a) Let $d=\lim_{c\to\infty}[\alpha_1(c)+\alpha_2(c)]<1$.
For any $F\in\mathbb{B}(\NN)$ define
$${\cal G}\circ F(c)=\min\{\alpha_1(c)F(c+1)+\alpha_2(c)F(c-1)+\alpha_3(c);\lambda\}.$$
Then, for any functions $F_1,F_2\in\mathbb{B}(\NN)$, it is easy to show that
\begin{align}
{\cal G}\circ F_1(c) & \le  {\cal G}\circ F_2(c)+ \max\{\alpha_1(c)[F_1(c+1)-F_2(c+1)]+\alpha_2(c)[F_1(c-1)-F_2(c-1)];0\}. \label{e36}
\end{align}
Therefore,
  $$|{\cal G}\circ F_1(c)-{\cal G}\circ F_2(c)|\le d \cdot\sup_{c\in{\NN}}|F_1(c)-F_2(c)|,$$
where $d=\sup_{c\in \NN} \{ \alpha_1(c)+\alpha_2(c)\}<1$ by hypothesis.
Consequently, ${\cal G}$ is a contraction in $\mathbb{B}(\NN)$, and equation (\ref{e34}), equivalent to $w={\cal G}\circ w$, has a single bounded solution which can be built by successive approximations
\begin{equation}\label{e37}
w_0(c)=0;~~~~w_{n+1}={\cal G}\circ w_n.
\end{equation}

Note that, since the functions $\alpha_j(c)$, $j=1,2,3$, are non-negative and increase with  $c\in\NN$ and, since $\lambda>0$, all the functions $w_n(c)$ increase with  $c\in\NN$ and $n\in\NN$.

(b)  Suppose $\lambda<\lambda^*$ and assume that $\forall c\in\NN$ the first expression in (\ref{e34}) is not bigger than $\lambda$. Then
  $$w(c)=\alpha_1(c)w(c+1)+\alpha_2(c) w(c-1)+\alpha_3(c),$$
so that $\lim_{c\to\infty}w(c)=\lambda^*>\lambda$, contradiction. Therefore, there is (the smallest) $c^*\in\NN$  such that the first expression in (\ref{e34}), which increases with $c$, is strictly bigger than $\lambda$ if and only if $c\ge c^*$. Note that $c^*>0$ because $\alpha_1(0)=\alpha_2(0)=\alpha_3(0)=0$.

Suppose $\lambda\ge\lambda^*$. Then one can easily show by induction that $w_n(c)\le \lambda^*$ for all $c\in\NN$ and $n\in\NN$. Hence the first expression in (\ref{e34}) does not exceed
  $$\frac{\zeta_3}{1-\zeta_1-\zeta_2}\cdot(\zeta_1+\zeta_2)+\zeta_3=\lambda^*\le\lambda.$$
$~$\hfill $\Box$

Let us go back to the epidemic model and make the following hypothesis.
\begin{hypot}
\item\label{HD} The functions
$$\alpha_1(c)=\frac{\rho_b(c)}{\eta+\rho_b(c)+\rho_d(c)+\kappa_i(c)},~~\alpha_2(c)=\frac{\rho_d(c)}{\eta+\rho_b(c)+\rho_d(c)+\kappa_i(c)}$$
  $$\mbox{ and }
\alpha_3(c)=\frac{\frac{\kappa_i(c)}{(\eta+\kappa_r)}}{\eta+\rho_b(c)+\rho_d(c)+\kappa_i(c)}$$
increase with $c$.
\end{hypot}

According to Lemma \ref{exl2}(a), equation (\ref{e33}) has a single bounded solution $v(c)$. From Lemma \ref{exl1} and Lemma \ref{exl2}(b) we deduce that $V(s,c,i)=sv(c)+\frac{i}{\eta+\kappa_r}$ and the optimal control strategy looks as follows. If
$$\lambda<\lambda^*=\lim_{c\to\infty}\frac{\alpha_3(c)}{1-\alpha_1(c)-\alpha_2(c)}=\lim_{c\to\infty}\frac{\frac{\kappa_i(c)}{(\eta+\kappa_r)}}{\eta+\kappa_i(c)},$$
then there exists $c^*\in\NN^*$  such that the first expression in (\ref{e32}) is strictly positive if and only if $c\ge c^*$. Thus
  $${\bf X}^g=\{(s,c,i)\in{\bf X}:~~c<c^*\},~~{\bf X}^i=\{(s,c,i)\in{\bf X}:~~c\ge c^*\}.$$
This means that one should immunize simultaneously all the existing susceptibles as soon as the number of carriers is bigger or equal $c^*$. This immunization is necessary also at the initial time moment zero in case $c_0\ge c^*$. If the number of carriers is below $c^*$ then wait.

If $\lambda\ge \lambda^*$ then $c^*=\infty$ and ${\bf X}^i=\emptyset$: one should never immunize any susceptibles since the cost of immunizations is too high.

It is interesting to note that  the critical value  $\lambda^*$ does not depend on the birth and death parameters $\rho_b(c)$ and $\rho_d(c)$. The optimal control strategy is of threshold type. Note also that all the susceptibles interact with the carriers independently of each other. Hence one can make a decision about the immunization of any one susceptible independently of the number of susceptibles, i.e. all of them are immunized (if needed) simultaneously. The cost coming from the current number of the infectives is not under control, hence the number of infectives does affect the optimal control strategy.

\bibliography{af}

\end{document}